\documentclass{amsart}
\pdfoutput=1

\usepackage{amsmath,amssymb}
\usepackage{lmodern,textcomp}  
\usepackage{esint}
\usepackage{enumerate}   
\usepackage{upgreek}
\usepackage{bbm}
\usepackage{mathrsfs}
\usepackage{cite}
\usepackage{xparse}
\usepackage[normalem]{ulem}

\usepackage{xcolor}
\usepackage{hyperref}
\usepackage[a4paper,margin=1.1in]{geometry}

\newtheorem{thm}{Theorem}[section]
\newtheorem{corollary}[thm]{Corollary}
\newtheorem{lemma}[thm]{Lemma}
\newtheorem{prop}[thm]{Proposition}
\theoremstyle{definition}

\theoremstyle{remark}
\newtheorem{rem}[thm]{Remark}
\numberwithin{equation}{section}
\newtheorem{example}[thm]{Example}

\newcommand{\R}{{\mathbb{R}}}
\newcommand{\ern}{{\mathbb{R}^N}}

\newcommand{\N}{\mathbb{N}}
\newcommand{\ep}{\varepsilon}
\newcommand{\ud}{\mathrm{d}}
\newcommand{\conv}{{\mathcal K}}

\newcommand{\Om}{{\Omega}}
\newcommand{\adm}{\mathscr{A}_{G}}

\newcommand{\Cc}{{\mathcal C}}

\newcommand{\M}{\mathscr{M}}

\newcommand{\Per}{\mathcal{P}}
\newcommand{\I}{\mathcal{I}}

\newcommand{\weakstar}{\overset{\ast}{\rightharpoonup}}
\newcommand{\TV}{\mathcal{TV}}

\DeclareMathOperator*{\dist}{{\rm dist}}

   \def\XXint#1#2#3{{\setbox0=\hbox{$#1{#2#3}{\int}$}
        \vcenter{\hbox{$#2#3$}}\kern-.5\wd0}}

\def\dy{\,{\rm d}y}
\def\dx{\,{\rm d}x}

\def\vs{\vspace{1mm}}

\renewcommand{\Gamma}{\varGamma}

\title[Asymptotics of fractional perimeters and first-order selection]
{{Asymptotics of fractional perimeters \\ and first-order selection of minimizers \\ with prescribed exterior datum}}

\author{L. De Luca}
\address[Lucia De Luca]{Istituto per le Applicazioni del Calcolo ``M. Picone'', IAC-CNR, Roma (Italy)}
\email{\url{lucia.deluca@cnr.it}}

\author{G. Palatucci}  
\address[Giampiero Palatucci]{Dipartimento di Scienze Matematiche, Fisiche e Informatiche, Universit\`a di Parma, Campus - Parco Area delle Scienze 53/A, 43124 Parma  (Italy)}  \email{\url{giampiero.palatucci@unipr.it}}

\author[M. Piccinini]{Mirco Piccinini} 
 \address[M. Piccinini]{Dipartimento di Matematica, Politecnico di Milano, Via Edoardo Bonardi~9, 20133, Milano (Italy)}
 \email{\url{mirco.piccinini@polimi.it}}

\author{M. Ponsiglione}
\address[Marcello Ponsiglione]{Dipartimento di Matematica ``G. Castelnuovo'', Sapienza Universit\`a di Roma, Roma (Italy)} 
\email{\url{ponsigli@mat.uniroma1.it}}

\begin{document}



\thanks{{\it Funding.} L.~De Luca has been supported by the INdAM project 2026 ``Ottimizzazione Geometrica e Frontiere Libere'' -- CUP\_E53C25002010001. G.~Palatucci has been supported by the University of Parma via FIL 2023 Project ``Nonuniform ellipticity in nonlinear problems'' CUP\_D93C24000380005.
G.~Palatucci and M.~Piccinini have been supported by the INdAM project 2026
``Parabolic equations with drift: well posedness, gradient estimates and
boundary effects'' -- CUP\_E53C25002010001.}

\begin{abstract}
We investigate $s$-perimeter functionals with prescribed exterior data in the singular limit~$s \to 0^+$. In this regime, although the leading-order asymptotics depends on the global distribution of the exterior datum, it may collapse to a purely volumetric functional and therefore fail to distinguish among different competitors. To capture the first nontrivial variational contribution beyond this degenerate scenario, 
we establish a full  $\Gamma$-convergence expansion; 
the resulting limit functional consists of a spatial interaction term measuring the residual interaction between the phase inside the domain, and a scalar far-field correction accounting for the prescribed exterior datum.

Applications to unconstrained and volume-constrained minimization problems show that perturbations, invisible at the zero-order, may nevertheless select the limiting shape of minimizers. This provides a selection principle for small-order fractional perimeter problems and makes explicit how the prescribed
exterior datum affects the variational problem as~$s\to0^+$.
\end{abstract}

\maketitle

{\tableofcontents}

\section{Introduction}\label{sec:intro}

Fractional perimeters model interfaces whose energy is not concentrated on
the boundary itself, but is produced by interactions across all length scales.

For $\Om\subset\R^N$  {  bounded open set with Lipschitz boundary}, $E\subset\R^N$  measurable,
and $s\in(0,1)$, the $s$-fractional perimeter of $E$ in $\Om$ is
\begin{equation*}\label{def_per_s}
{\rm Per}_s(E,\Om) := { \int_{E \cap \Om}\ud x\int_{\Cc E}\frac{\ud y}{|x-y|^{N+s}} +\int_{\Cc E \cap \Om}\ud x\int_{E \cap \Cc \Om}\frac{\ud y}{|x-y|^{N+s}}},
\end{equation*}
where { for~$A \subset \R^N$ we let~$\Cc A := \R^N \setminus A$.}

This functional was introduced by Caffarelli, Roquejoffre and Savin
in~\cite{CRS10} as a nonlocal counterpart of the classical perimeter. Since
then, nonlocal minimal surfaces have developed into a rich
geometric-variational theory. The available results include uniform estimates
and limiting arguments~\cite{CV11,CV12}, bootstrap regularity methods~\cite{BFV14},
classification of low-dimensional cones~\cite{SV13}, Bernstein-type results and gradient estimates for nonlocal
minimal graphs~\cite{FV17,CC19}, stability and quantitative flatness~\cite{CSV19},
and constructions of nontrivial stationary or minimizing examples, such as
nonlocal Delaunay surfaces and Lawson cones~\cite{DDPV16,DPW18}.
For $\Omega=\R^N$, fractional perimeters appeared in connection with long-range dissipative phenomena in~\cite{V90,V91}, in nonlocal curvature flows~\cite{I09,CS10,CMP15,CDNP21}, while isoperimetric and stability properties have been studied in~\cite{FMM11,FFMMM15}.
More broadly, geometric and free-boundary problems driven by nonlocal integro-differential operators have given rise to parallel regularity theories~\cite{CROS17}.

The dependence on the fractional parameter~$s$ is especially delicate.
Near the local regime, say~$s\sim1$, flatness results are known in
dimension three~\cite{CCS20}. More recently, the Simons cone in
$\R^8$ has been proved to minimize the fractional perimeter under
compact perturbations for all $s$ sufficiently close to~$1$~\cite{NP26}.
At the highly nonlocal end, namely~$s\sim0$, one finds genuinely
different phenomena, including stable nonlocal Lawson cones in
dimension seven~\cite{DPW18} and the recent rigidity of stable
$s$-minimal cones in the plane~\cite{Cas25}.
In particular, a Lawson cone in $\R^7$ associated with the $3+4$
splitting has recently been proved to minimize the fractional
perimeter for all sufficiently small positive orders~$s$~\cite{Pal26}.

\medskip
Accordingly, the two limit regimes of the perimeter functional as $s\to1^-$ and~$s\to0^+$ are fundamentally different.
At the local endpoint,
after multiplication by~$(1-s)$, the fractional
perimeter recovers its classical counterpart, up to a dimensional
factor. This is the content of~\cite{ADM11} where it has been shown that
\[
 (1-s){\rm Per}_s(\cdot,\Om)
\overset{\Gamma}{\longrightarrow}
\omega_{N-1}\mathcal H^{N-1}(\partial^\ast \cdot\cap\Om).
\]

In Section~\ref{sec:Gamma1} we revisit such analysis, fixing the exterior datum instead of perturbing it using boundary layers killing the contributions on~$\partial\Omega$. The corresponding $\Gamma$-convergence result is Theorem~\ref{s1}, which states 
\begin{equation*}\label{intro_s-a-1}
 (1-s){\rm Per}_s(\cdot,\Om)
\overset{\Gamma}{\longrightarrow}
\omega_{N-1}\mathcal H^{N-1}(\partial^\ast \cdot\cap\overline\Om).
\end{equation*}

\medskip
The opposite endpoint $s\to0^+$ is not a regular perturbation of the local
perimeter problem. The kernel approaches the borderline power~$|x-y|^{-N}$,
and the far-field integral diverges at the scale~$1/s$. Thus the same
functional carries two difficulties at once: a singular near-field interaction
and a degenerate far-field contribution. In the presence of a prescribed
exterior datum, the far field is not a harmless additive constant. It competes
with the phase inside~$\Om$, and its leading effect is governed by a global tail (volume) balance
rather than by the geometry of the interface. The first step is therefore to
identify this collapse; the central variational question is what survives
after it is subtracted.

Far-field interactions are a central feature of nonlocal analysis. In fractional
equations they enter the regularity theory through the nonlocal tails
introduced in~\cite{DCKP14,DCKP16}. For fractional perimeters, the effective
tail at infinity was isolated in~\cite{DFPV13} in the  pointwise limit of
the $s$-perimeter as $s\to0^+$. The same endpoint is governed by global,
rather than purely local, effects also in the asymptotic analysis of
fractional perimeters on Riemannian manifolds~\cite{CG24}. A recent
complementary perspective is given in~\cite{DFP26}, where the
Maz'ya--Shaposhnikova asymptotics~\cite{MS02} are extended to localized
Gagliardo seminorms for functions with nontrivial mass at infinity, recovering
in particular the pointwise asymptotics of fractional perimeters on
characteristic functions.

\medskip
The present paper takes this leading-order picture as its starting point and
studies the first variational layer at which the geometry lost in the
tail-volume limit reappears. At this scale, the limiting functional yields a
genuine selection principle for minimizers in regimes where the leading-order functional is degenerate.

 The regime considered here is not a localization of the whole-space critical problem. In the latter, the energy is global and symmetric, and the critical
object is governed essentially by the interaction between a set of finite
volume and its complement.
 This is the framework of the~$0$-fractional perimeter studied 
in~\cite{DLNP21} (see also~\cite{CDKNP23}), where Riesz-type functionals naturally appear. 
Related critical structures also arise in the analysis of the logarithmic Laplacian,
introduced in~\cite{CW19} as the derivative at~$s=0$ of the fractional
Laplacian. In the same endpoint regime,
Jarohs, Sen and Weth~\cite{JSW26} analyzed $s$-harmonic functions with prescribed exterior data and showed that both the zero-order asymptotics and the differentiability with respect to~$s$ of the associated Poisson problem can be expressed in terms of the logarithmic Laplacian of suitable extensions of the datum.

Here, by contrast, the phase is fixed outside the domain, the exterior datum
may have infinite measure, and the order-one expansion  must distinguish
internal interactions, interactions with the prescribed datum, interactions
with the exterior complement, and boundary terms generated by the cut at
$\partial\Om$. These contributions do not separate automatically. The
renormalization must remove the finite-$s$ divergent tail at the scale at
which it is produced, while leaving precisely the critical spatial interaction
associated with the kernel $|x-y|^{-N}$.

\medskip
Now we describe in more detail our results. Throughout the paper~$\Om\subset\R^N$ is a bounded open set with Lipschitz
boundary, and~$G\subset\Cc\Om$ is a measurable exterior datum.

We are interested in the following minimum problem,
$$
\min\bigl\{{\rm Per}_s(E,\Om):\ E\in\M(\R^N), \ E\cap\Cc\Om=G\bigr\}, 
$$
and more precisely in its behavior for $s\sim 0$.
To this end, we develop a $\Gamma$-convergence analysis for~${\rm Per}_s$.

\smallskip

The zero-order scaling is~$s{\rm Per}_s$. 
It kills the interactions
contained in bounded regions and retains only the exterior far field, producing
an affine functional of~$|E\cap\Om|$ whose coefficients are the two exterior
tail masses. This accounts at once for the nondegenerate and degenerate cases:
an unbalanced datum favours one pure phase, whereas balanced data, or any
fixed-volume class, leave the leading order constant.

The first-order normalization subtracts the finite-$s$ tail-volume term itself,
before taking the limit. This cancellation removes the exterior divergence at
its natural scale and exposes a critical spatial interaction with
kernel~$|x-y|^{-N}$. Since the phase is prescribed outside~$\Om$, the resulting
functional is not a whole-space critical perimeter restricted to~$\Om$: the cut
at~$\partial\Om$ and the exterior datum remain structural parts of the limit.

\medskip
Let us also indicate the main points of the argument. At order zero, the
far-field contribution is isolated by comparing~$|x-y|^{-N-s}$ with~$|y|^{-N-s}$ for large~$|y|$. At order one, after the tail cancellation, the
renormalized energies provide localized nonlocal control of~$\chi_{E\cap\Om}$; combined with the compact embedding theorem of Jarohs and Weth~\cite{JW20}, this yields strong~$L^1$ compactness. The liminf and recovery
arguments are then carried out by decomposing the energy into internal,
near-exterior and far-exterior contributions, keeping the exterior datum fixed
throughout.

\medskip
\subsection{Main results and organization of the paper}\label{sez_mains}

We now describe the main results in the order in which they enter the
analysis: zero- and first-order~$\Gamma$-convergence, stability with respect to the exterior datum,
minimization, and comparison with the local endpoint.

\medskip
\noindent{\it Zero-order collapse.}
Assume, when needed, that the following limit
$$
\alpha(G):=\lim_{s\to0^+}\alpha_1^s(G) = \lim_{s\to0^+}s\int_{G\setminus B_1}\frac{\ud z}{|z|^{N+s}}\,,
$$
exists.
 Since $\alpha(\Cc G)=N\omega_N-\alpha(G)$, the natural tail-volume
functional is
$$
\mu_G(m):= \alpha(\Cc G)m+\alpha(G)(|\Om|-m), \qquad m\in[0,|\Om|].
$$

The first result identifies the leading-order asymptotics. For every
vanishing sequence $\{s_k\}_{k\in\N}\subset(0,1)$, we prove that
\begin{equation*}\label{intro_gamma0}
s_k\Per^{s_k}_G(\cdot)
\overset{\Gamma}{\longrightarrow}
\mu_G\bigl(\|\cdot\|_{L^1(\Om)}\bigr)\,,
\end{equation*}
with respect to weak-$\ast$ convergence of measures
$$
\chi_{E_k\cap\Om}\,\ud x\weakstar f\,\ud x, \qquad 0\le f\le1.
$$
This scalar limit is Theorem~\ref{thm:G-lim} in
Section~\ref{sec:gamma0th}.

\medskip
\noindent{\it First-order $\Gamma$-limit.}
The first-order limit is the functional
$$
\Per^0_G(E)=\widetilde\Per^0(E)+\widehat\Per^0_{1,G}(E).
$$
The term $\widetilde\Per^0$ (defined in \eqref{tilde}) is the internal critical interaction in
$\Om$, while $\widehat\Per^0_{1,G}$ (defined in \eqref{cappuccio}) contains the interaction with the
exterior datum together with the renormalized tail contribution. Thus
$\Per^0_G$ retains the boundary cut and the datum as structural parts of the
limit, rather than compressing them into a scalar tail mass.

For the finite-$s$ tail-volume contribution we set
$$
\mu^{(s)}_G(m) := \alpha_1^s(\Cc G)m+\alpha_1^s(G)(|\Om|-m),
\qquad m\in[0,|\Om|].
$$

The main first-order $\Gamma$-convergence result states that, for every
vanishing sequence $\{s_k\}_{k\in\N}\subset(0,1)$,
$$
 \Per^{s_k}_G(\cdot) - \frac1{s_k}\mu^{(s_k)}_G(|\cdot\cap\Om|)
\overset{\Gamma}{\longrightarrow}
 \Per^0_G(\cdot)\,,
$$
with respect to strong $L^1(\Om)$ convergence. This is Theorem~\ref{thmprimoord}. 

The use of the finite-$s$ quantity~$\mu^{(s)}_G$ is essential: subtracting only its limit would not remove the
tail contribution at the scale where it is actually produced.

When the scalar correction
\begin{equation*}\label{intro_def_lambda}
\lambda_G:= \lim_{s\to0^+} \frac{\alpha_1^s(G)-\alpha(G)}{s},
\end{equation*}
exists, the same first-order result can be written as an expansion around the
zero-order tail mass:
\begin{equation*}\label{intro_gamma1st}
\Per^{s_k}_G(\cdot) - \frac1{s_k}\mu_G(|\cdot\cap\Om|) \overset{\Gamma}{\longrightarrow}
 \mathcal F_G^{(1)}(\cdot),
\end{equation*}
where
\begin{equation*}\label{intro_limite-1st}
\mathcal F_G^{(1)}(E) :=\Per^0_G(E) +\lambda_G\bigl(|\Cc E\cap\Om|-|E\cap\Om|\bigr).
\end{equation*}
This is Theorem~\ref{teolam}. The formula separates the first-order tail bias
from the spatial part of the limiting energy.

\medskip
\noindent{\it Volume constraints and stability under varying exterior data.}
The same analysis is carried out under the constraint $|E\cap\Om|=m$. In
this case the finite-$s$ tail term $\mu^{(s)}_G(m)$ is constant on the
admissible class, and the limit is the restriction of $\Per^0_G$ to
competitors with fixed volume. The corresponding compactness and
$\Gamma$-convergence result is Theorem~\ref{teovol0}.

We also prove that the first-order theory is stable under variations of the
exterior datum. If~$G_k\subset\Cc\Om$ converges to~$G$ in the sense that
$|G_k\triangle G|\to0$,
 then the renormalized functionals with exterior data~$G_k$ converge to
the first-order functional associated with~$G$. The combined variant accounting both for volume-constraints and varying exterior data is included as well. These are Theorems~\ref{teocambG} and~\ref{teoval0var}.

\medskip
\noindent{\it Selection of minimizers.}
The final part of the paper applies the theory to fractional isoperimetric
problems with prescribed exterior data. In the unconstrained problem, the
leading tail mass selects one of the two pure phases, except in the balanced
case
$$
\alpha(G)=\frac{N\omega_N}{2},
$$
where the zero-order limit is completely degenerate and the first-order
functional gives the remaining selection principle. Under a prescribed volume
constraint, both the leading term $\mu^{(s)}_G(m)$ and the scalar correction
$\lambda_G(|\Om|-2m)$ are constant. Hence shape selection can only be
carried by the spatial part of $\Per^0_G$. The fixed-volume consequences of
the $\Gamma$-convergence are collected in Proposition~\ref{prop_volume-min}.

The selection mechanism is then made explicit through balanced exterior data
with bounded perturbations. Starting from a balanced cone $G_A$ (see  Formula~\eqref{def_cono_a}), we add a
bounded window $P$ and remove a bounded hole $Q$. These perturbations do not
change the leading tail balance,
but they enter the first-order functional through the spatial potential
$$
\Phi_{P,Q}(x) = \int_{\Cc B_1}
(\chi_P(y)-\chi_Q(y))\left(\frac1{|x-y|^N}-\frac1{|y|^N}\right)\ud y.
$$
For $F=E\cap\Om$, the corresponding decomposition is
$$
\Per^0_G(E) = { \Per^0_{G_A}(F\cup G_A)} + \int_\Om\Phi_{P,Q}\,\ud x - 2\int_F\Phi_{P,Q}\,\ud x.
$$
Consequently, bounded exterior windows and holes can be invisible to the
leading tail mass and, in a fixed-volume class, can contribute only a constant
scalar bias, while still selecting different limiting minimizers through their
first-order spatial tail geometry. This is proved in
Propositions~\ref{prop_win-hol_bilanciate}.

\medskip
\noindent{\it Comparison with the local endpoint.}
We finally compare this selection mechanism with the local regime
$s\to1^-$. If the windows and holes stay away from~$\Om$, then they are
invisible in the local limit, since Theorem~\ref{s1} only sees the exterior
trace near $\partial\Om$. The same perturbations can instead be decisive
as $s\to0^+$, because the vanishing-order regime is governed by the
renormalized tail. In the examples with moving exterior data, diagonal
regimes near $s=1$ are controlled by the competition between the local scale
and the distance of the datum from~$\partial\Om$, whereas the regime~$s\to0^+$ is controlled by the tail interaction.

\medskip
 Thus the two endpoint regimes have different variational sensitivities:
exterior data that are invisible in the local limit may determine the
first-order selection mechanism as~$s\to0^+$.

%
%

\medskip
\section{Preliminaries}\label{sec:preli}
Let $\M(\R^N)$ denote the family of Lebesgue measurable subsets of $\R^N$.
For any set~$A \in \M(\R^N)$ we denote by $\Cc A$ the complement set of~$A$ in~$\R^N$. For every
$s\in[0,1)$ and every $A,B\in\M(\R^N)$ we set
\begin{equation}\label{intera}
\I^s(A,B):=\int_A\ud x\int_B \frac{\ud y}{|x-y|^{N+s}}.
\end{equation}
The value $+\infty$ is allowed.

Throughout the paper, $\Om\subset\R^N$ is a bounded open set with Lipschitz
boundary. Given $s\in(0,1)$, the $s$-fractional perimeter of a measurable set
$E\in \M(\R^N)$ relative to $\Om$ is
\begin{equation}\label{per}
{\rm Per}_s(E,\Om)
:=\I^s(E\cap\Om,\Cc E) + \I^s((\Cc E)\cap\Om,E\cap\Cc\Om).
\end{equation}
We refer to Chapters~V--VI in~\cite{FAG} for an account of the fractional perimeter and
its exterior data.
\smallskip

Let $G\in\M(\R^N)$ be such that $G\subset\Cc\Om$. Throughout the paper,~$G$ represents the prescribed exterior datum; we introduce the admissible class 
\begin{equation*}\label{def:admissible}
\adm :=\Big\{ E\in\M(\R^N):E\cap\Cc\Om=G\Big\}.
\end{equation*}
For every $R>0$ and every $s\in(0,1)$ we define
\begin{equation}\label{alfaiso}
\alpha^s_R(G) = sR^{-s}\,\operatorname{Tail}_{s/2,2}(\chi_G;0,R) = s\int_{G\cap\Cc B_R}\frac{\ud z}{|z|^{N+s}},
\end{equation}
where $\operatorname{Tail}_{\sigma,p}$ is the nonlocal tail introduced in \cite{DCKP14,DCKP16}.
Here and below, for every $\rho>0$ and for every $\xi\in\R^N$, we denote by $B_\rho(\xi)$ the open ball with radius $\rho$ and center $\xi$ and we set $B_\rho:=B_\rho(0)$.

Following~\cite{DFPV13}, for the zero-order analysis in Section~\ref{sec:gamma0th}, we assume that the leading tail mass {$\alpha (G)$} exists, where 
\begin{equation}\label{eq:alphaG}
\alpha(G):=\lim_{s\to 0^+}\alpha_1^s(G).
\end{equation}
This quantity measures the leading-order contribution of the exterior datum
at infinity in the limit $s\to0^+$. Notice that
$$
0\le\alpha(G)\le N\omega_N,
$$
where $\omega_N$ denotes the volume of the unit ball in $\R^N$.

\begin{rem}\label{usoord0}
The quantities $\alpha^s_R(G)$ and~$\alpha(G)$ are actually defined for any measurable
set $G\in \M(\R^N)$; the condition $G\cap\Om=\emptyset$ is only relevant for
the variational problem with prescribed exterior datum.

As highlighted in~\cite[Observation 2]{DFPV13}, if the limit $\alpha(G)$
exists, then
\begin{equation}\label{eq:alphaG2}
\alpha(G)=\lim_{s\to 0^+}\alpha_R^s(G),
\qquad\textrm{for every }R>0.
\end{equation}
Thus, at leading order, the reference radius is immaterial.

Finite volume perturbations do not affect the value of~$\alpha(G)$.
 Indeed, if~$G,G'\in\M(\R^N)$ and~$|G\triangle G'|<+\infty$, then, for every~$s\in(0,1)$ and every~$R>0$,
\begin{equation}\label{pertufin}
\big|\alpha^s_R(G)-\alpha^s_R(G')\big|
\le s\int_{(G\triangle G')\cap\Cc B_R}\frac{\ud z}{|z|^{N+s}}
\le sR^{-N-s}|G\triangle G'|.
\end{equation}
Consequently, whenever $\alpha(G)$ exists, also $\alpha(G')$ exists and
\begin{equation}\label{pertufin2}
\alpha(G')=\alpha(G).
\end{equation}
In particular, $\alpha(G)=0$ whenever $G$ has finite volume.

Finally, since for every $s\in(0,1)$ and every $R>0$,
\begin{equation}\label{alpharcomp}
\alpha^s_R(\Cc G)=N\omega_NR^{-s}-\alpha^s_R(G),
\end{equation}
if~$\alpha(G)$ exists, then also $\alpha(\Cc G)$ exists and
\begin{equation}\label{eq:CG-alpha}
\alpha(\Cc G)=N\omega_N-\alpha(G).
\end{equation}
 We say that $G$ is {\it balanced} if 
\begin{equation}\label{balcase}
\alpha(G)=\frac{N\omega_N}{2}.
\end{equation} 
In such a case, $\alpha(\Cc G)=N\omega_N/2$ as well,
so that  the balanced condition 
is invariant
under complementation.

\end{rem}

Whenever the limit $\alpha(G)$ in \eqref{eq:alphaG} exists, we may introduce the following quantity, which measures in some sense the rate of convergence of $\alpha^s_1(G)$ to $\alpha(G)$,
\begin{equation}\label{lamG}
\lambda_G := \lim_{s\to0^+}\frac{\alpha^s_1(G)-\alpha(G)}{s}.
\end{equation}
When $\lambda_G\in\R$, we refer to it as  the scalar first-order tail correction associated 
with the exterior datum $G$. In contrast with $\alpha(G)$, one may easily see that compact perturbations may affect the scalar correction $\lambda_G$ so that its definition is sensitive to the choice of the radius; in other words, $\lambda_G$ may vary when replacing $\alpha^s_R(G)$ with $\alpha^s_1(G)$ in \eqref{lamG}.
\medskip

%
%

\medskip

\section{The zero-order $\Gamma$-limit} \label{sec:gamma0th}
This section is devoted to the $\Gamma$-convergence of the functionals
{$s{\rm Per}_s(\cdot,\Omega)$}
under fixed boundary datum $G$.

To this end, for every $s\in(0,1)$ we define ${\Per^s_{G}}:\M(\ern)\to[0,+\infty] $ as
\begin{equation}\label{nuovoper}
{\Per}^s_{G}(E):=\left\{\begin{array}{ll}
{\rm Per}_s(E,\Om)&\,\textrm{if }E\in\adm,\\
+\infty&\,\textrm{otherwise}.
\end{array}
\right.
\end{equation}
{Therefore, for every $E\in\adm$, we have
\begin{equation}\label{nuovoper1}
{\Per}^s_{G}(E)=\I^s(E\cap\Omega,\Cc E\cap\Omega)+\I^s(E\cap\Omega,\Cc G\cap\Cc\Omega)+\I^s(\Cc E\cap\Omega,G).
\end{equation}
Throughout this section, we assume that the limit $\alpha(G)$ in \eqref{eq:alphaG} exists.} 
At leading order, the nonlocal perimeter only detects how much of $\Om$ is occupied by the phase and how much is left to its complement, as already seen in the analysis in~\cite{DFPV13}. This motivates the following affine $\Gamma$-limit candidate $\mu_G:[0,|\Om|]\to[0,+\infty)$:
\begin{equation}\label{gammaorzero}
\mu_{G}(m):=\alpha(\Cc G)m+\alpha(G)(|\Om|-m).
\end{equation}
Note that  in the balanced case \eqref{balcase}
\begin{equation}\label{costa}
\mu_{G}(m)\equiv \frac{N\omega_N}{2}|\Om|\,,\qquad\textrm{for every }m\in[0,|\Om|].
\end{equation}
As a consequence, the leading-order tail mass does not select between the two phases in $\Om$. 

Our first main result identifies the zero-order asymptotics of $s{\Per}^s_{G}$.
\begin{thm}\label{thm:G-lim}
For every $s\in(0,1)$, let ${\Per}^s_{G}$ and $\mu_G$ be defined in \eqref{nuovoper} and \eqref{gammaorzero}, respectively.
Let $\{s_k\}_{k\in\N}\subset(0,1)$ be a vanishing sequence.
The following $\Gamma$-convergence result holds true.
\begin{itemize}
\item[(i)] ($\Gamma$-liminf inequality) Let $f\in  L^{\infty}(\Om;[0,1])$ and let $\{E_k\}_{k\in\N}\subset\M(\ern)$ be such that $\chi_{E_k\cap\Om}\ud x\weakstar f\ud x$. Then
\begin{equation}\label{liminf}
\mu_{G}(\|f\|_{L^1})\le\liminf_{k\to +\infty}s_k\Per^{s_k}_{G}(E_k).
\end{equation}
\item[(ii)] ($\Gamma$-limsup inequality) For every $f\in  L^{\infty}(\Om;[0,1])$ there exists $\{E_k\}_{k\in\N}\subset\adm$ with $\chi_{E_k\cap\Om}\ud x\weakstar f\ud x$ such that
\begin{equation*}\label{limsup0}
\mu_{G}(\|f\|_{L^1})=\lim_{k\to +\infty}s_k\Per^{s_k}_{G}(E_k).
\end{equation*}
\end{itemize}
\end{thm}
\begin{rem}\label{trivia}
\rm{
Notice that given a sequence $\{f_k\}_{k\in\N}\subset L^\infty(\Om;[0,1])$, there exists a function $f\in L^\infty(\Om;[0,1])$ and a not-relabeled subsequence such that $f_k\ud x\weakstar f\ud x$.
Therefore, the $\Gamma$-convergence result in Theorem \ref{thm:G-lim} can be accompanied with a trivial compactness property.
}
\end{rem}
{Before providing the proof of Theorem \ref{thm:G-lim} we highlight some basic facts that will be useful for this scope and throughout all the paper.}
\begin{rem}\label{lipsets}
\rm{
Notice that for every bounded set $F$ with Lipschitz continuous boundary, for every $s\in[0,1)$ and for every {bounded open set $A$
\begin{equation*}
\I^{s}(F,\Cc F\cap A)<+\infty.
\end{equation*} 
Moreover,
$$
\lim_{s\to 0^+}\I^{s}(F,\Cc F\cap A)=\I^{0}(F,\Cc F\cap A).
$$
}
}
\end{rem}
\begin{rem}\label{contocruc}
\rm{{
Let $R>0$ be such that $\Omega\subset B_R$. For every $0<s<1$ and for every $E\in\adm$ we write
\begin{equation}\label{crucc1}
\begin{aligned}
&\,\I^s(E\cap\Omega, \Cc G\cap \Cc\Omega)-\frac{\alpha^s_R(\Cc G)}{s}|E\cap\Omega|\\
=&\,\I^s(E\cap\Omega,\Cc G\cap\Cc\Omega\cap B_R)+\int_{E\cap\Omega}\ud x\int_{ \Cc G\cap \Cc B_R}\bigg(\frac{1}{|x-y|^{N+s}}-\frac{1}{|y|^{N+s}}\bigg)\,\ud y
,
\end{aligned}
\end{equation} 
so that, by Remark~\ref{lipsets} above, we get
\begin{equation}\label{stimaalto1}
\begin{aligned}
&\,\bigg|\I^s(E\cap\Omega, \Cc G\cap \Cc\Omega)-\frac{\alpha^s_R(\Cc G)}{s}|E\cap\Omega|\bigg|\\
\le&\,\I^s(\Omega,\Cc\Omega\cap B_R)+\int_{\Omega}\ud x\int_{\Cc B_R}\bigg|\frac{1}{|x-y|^{N+s}}-\frac{1}{|y|^{N+s}}\bigg|\,\ud y
\\
\le&\, C(N,\Omega,R).
\end{aligned}
\end{equation}
Analogously,
\begin{equation}\label{crucc2}
\begin{aligned}
&\,\I^s(\Cc E\cap\Omega, G)-\frac{\alpha^s_R(G)}{s}|\Cc E\cap\Omega|\\
=&\,\I^s(\Cc E\cap\Omega,G\cap B_R)+\int_{\Cc E\cap\Omega}\ud x\int_{  G\cap \Cc B_R}\bigg(\frac{1}{|x-y|^{N+s}}-\frac{1}{|y|^{N+s}}\bigg)\,\ud y,
\end{aligned}
\end{equation}
and
\begin{equation}\label{stimaalto2}
\bigg|\I^s(\Cc E\cap\Omega, G)-\frac{\alpha^s_R(G)}{s}|\Cc E\cap\Omega|\bigg|
\le C(N,\Omega,R).
\end{equation}
}
}
\end{rem}

\begin{proof}[Proof of Theorem {\rm\ref{thm:G-lim}}]
We first prove the $\Gamma$-liminf inequality (i).
We can assume without loss of generality that
\begin{equation*}\label{enbou}
\sup_{k\in\N}s_k\Per^{s_k}_{G}(E_k)\le M,
\end{equation*}
for some $M>0$.
{
Let $R>0$ be such that $\Omega\subset B_R$.
By \eqref{stimaalto1}, we have 
\begin{equation}\label{mults1}
\lim_{k\to+\infty}\Big(s_k\I^{s_k}(E_k\cap\Omega, \Cc G\cap \Cc\Omega)-\alpha^{s_k}_R(\Cc G)|E_k\cap\Omega|\Big)=0,
\end{equation}
whence, by \eqref{eq:alphaG2} and by lower semicontinuity, we deduce
\begin{equation}\label{eq:lw1lim}
\liminf_{k\to+\infty}s_k\I^{s_k}(E_k\cap\Omega, \Cc G\cap \Cc\Omega)\ge \alpha(\Cc G)\int_{\Omega}f\,\ud x.
\end{equation}
Analogously, using \eqref{stimaalto2}, we get
\begin{equation}\label{mults2}
\lim_{k\to+\infty}\Big(s_k\I^{s_k}(\Cc E_k\cap\Omega,  G)-\alpha^{s_k}_R( G)|\Cc E_k\cap\Omega|\Big)=0,
\end{equation}
and hence
\begin{equation}\label{eq:lw2lim}
\liminf_{k\to+\infty}s_k\I^{s_k}(\Cc E_k\cap\Omega, G)\ge \alpha( G)\int_{\Omega}(1-f)\,\ud x.
\end{equation}
}

By \eqref{nuovoper1}, \eqref{eq:lw1lim} and \eqref{eq:lw2lim}, we get \eqref{liminf}.

We now pass to the proof of (ii).
If $f=\chi_{F}$ for some $F\in\M(\ern)$ with $F\subseteq\Om$ having Lipschitz continuous boundary, 
we claim that the recovery sequence is provided by $E_k:=F\cup G$ (for every $k\in\N$). 
{
Indeed, by Remark \ref{lipsets} (applied with $A=\Omega$), we have
\begin{equation}\label{dentro}
\lim_{k\to +\infty}s_k\I^{s_k}(E_k\cap\Omega,\Cc E_k\cap\Omega)=\lim_{k\to +\infty}s_k\I^{s_k}(F,\Cc F\cap\Omega)=0.
\end{equation}
Moreover, let $R>0$ be such that $\Om\subset B_{R}$.
By \eqref{mults1} and \eqref{mults2}, we have that
\begin{equation}\label{fuori}
\begin{aligned}
&\,\lim_{k\to+\infty}s_k\I^{s_k}(E_k\cap\Omega, \Cc G\cap \Cc\Omega)=\alpha(\Cc G)|F\cap\Omega|,
\\
&\,\lim_{k\to+\infty}s_k\I^{s_k}(\Cc E_k\cap\Omega, G)=\alpha(G)|\Cc F\cap\Omega|.
\end{aligned}
\end{equation}

By \eqref{nuovoper1}, summing \eqref{dentro} and \eqref{fuori}, we get the claim for $f$ as above.
}

\smallskip  
The general case $f\in L^\infty(\Om;[0,1])$ follows by a standard density argument in $\Gamma$-convergence: 
one first approximates $f$ in the weak$^*$ topology of measures by characteristic functions of subsets of $\Om$ with Lipschitz continuous boundary, then, applies the construction above to each approximating set. Since $\mu_{G}$ is affine and continuous with respect to the weak$^*$ convergence of measures, the conclusion follows by diagonalization.
\end{proof}
\begin{rem}\label{rem:stable}
We briefly record that Theorem \ref{thm:G-lim} is stable under perturbations of the boundary datum. 
More precisely, let $G\in\M(\ern)$ and $\{G_k\}_{k\in\N}\subset\M(\ern)$ be such that
\begin{equation*}\label{GnG}
\sup_{k\in\N}|G_k\triangle G|\le C.
\end{equation*}
Then, by arguing verbatim as in the proof of Theorem~\ref{thm:G-lim} and using Remark \ref{usoord0},  one obtains the corresponding~$\Gamma$-convergence statement with~$G_k$ in place of~$G$. 
\end{rem}

%
%

\medskip
\section{The first order $\Gamma$-limit} \label{sec:gamma1st}
 We now pass to the first-order regime.  More precisely, we subtract from the $s$-fractional perimeter $\Per^s_G$  the scaled
zero-order contribution associated with the prescribed tail datum $G$, that is  $\frac{1}{s} \mu^{(s)}_G$   defined in \eqref{mur} below. 
In this way the leading
order term, which only depends on the volume fraction inside $\Om$, is
removed, and the next nontrivial contribution can be detected.

It is useful to introduce this renormalization with an auxiliary radius
$R>0$. The reason is that the estimates below are first proved at a fixed
renormalization scale, where the tail contribution is separated from the
local part in a transparent way. The dependence on $R$ will then be shown to
be uniformly controlled, so that one can eventually return to the unit-radius
normalization.

To this end, for every $R>0$ and for every $s\in(0,1)$ we define
$\mu^{(s)}_{R,G}:[0,|\Om|]\to[0,+\infty)$ as
\begin{equation}\label{mur}
\mu^{(s)}_{R,G}(m):=\alpha^s_R(\Cc G)m+\alpha^s_R(G)(|\Om|-m),
\end{equation}
and we set $\mu^{(s)}_{G}(\cdot):=\mu^{(s)}_{1,G}(\cdot)$.

Note that, if the limit in \eqref{eq:alphaG} exists, then, for every $R>0$,
and for every $m\in[0,|\Om|]$ 
$$
\lim_{s\to 0^+}\mu^{(s)}_{R,G}(m)=\mu_{G}(m),
$$
where $\mu_{G}$ is defined in \eqref{gammaorzero}.

The next elementary observation shows that changing the renormalization
radius only produces an explicitly controlled error. Indeed, given
$s\in (0,1)$, for every $0<R_1<R_2$ we have
$$
\alpha^s_{R_1}(G)=\alpha_{R_2}^{s}(G)+s\int_{G\cap (B_{R_2}\setminus B_{R_1})}\frac{\ud z}{|z|^{N+s}},
$$
and
$$
\alpha^s_{R_1}(\Cc G)=\alpha_{R_2}^{s}(\Cc G)+s\int_{\Cc G\cap (B_{R_2}\setminus B_{R_1})}\frac{\ud z}{|z|^{N+s}};
$$
it follows that
for every $m\in[0,|\Om|]$ 
 \begin{eqnarray}\label{usoora}
\phantom{+}\mu^{(s)}_{R_1,G}(m)-\mu^{(s)}_{R_2,G}(m)
& = & \phantom{+}sm\int_{\Cc G\cap (B_{R_2}\setminus B_{R_1})}\frac{\ud z}{|z|^{N+s}}\notag\\
&&+s(|\Om|-m)\int_{G\cap (B_{R_2}\setminus B_{R_1})}\frac{\ud z}{|z|^{N+s}}\notag\\
& \le & \phantom{+}s|\Om|\int_{B_{R_2}\setminus B_{R_1}}\frac{\ud z}{|z|^{N+s}} \ = 
{s|\Om|N\omega_N\frac{R_1^{-s}-R_2^{-s}}{s}}.
\end{eqnarray}
\\
As a consequence, for every $R>1$, for $s$ small enough we have
\begin{equation}\label{usoora2}
|\mu^{(s)}_{G}(m)-\mu^{(s)}_{R,G}(m)|\le 2sN\omega_N\log R\,|\Om|.
\end{equation}

\medskip
\subsection{Compactness}
 In this section we prove the compactness statement of the first-order
$\Gamma$-convergence analysis. The main point is that, after subtracting the zero-order
renormalized volume contribution, a uniform bound on the energy still controls
a genuinely compact nonlocal seminorm at the critical exponent.  As a consequence, in contrast with the zero-order analysis, any limit of  characteristic functions with uniformly bounded renormalized $s$-fractional perimeters is still  a characteristic function.

The proof relies on a compact embedding theorem for suitable nonlocal quadratic forms. We recall the following result from~\cite{JW20} (see also~\cite{CDKNP23,DPV12}).

\begin{thm}[Theorem~1.2 in~\cite{JW20}]\label{JaWe}
Let $A$ be a bounded open subset of $\ern$ with Lipschitz continuous boundary.
Let $\kappa:\ern\to[0,+\infty]$ be a radially symmetric kernel satisfying
\begin{equation}\label{iposulkernel}
\int_{\ern}\kappa(z)\,\ud z=+\infty\qquad\textrm{and}\qquad
\int_{\ern}\min\{1,|z|^2\}\kappa(z)\,\ud z<+\infty,
\end{equation}
and let
\begin{multline*}
\mathscr{W}^{\kappa}(A):=\Big\{f\in L^2(\ern)\,:\,f\equiv 0\textrm{ in }\Cc A,
\iint_{\ern\times\ern}|f(x)-f(y)|^2\kappa(x-y)\,\ud y\,\ud x<+\infty\Big\}
\end{multline*}
be the Banach space endowed with the norm

\begin{equation*}
\|f\|_{\mathscr{W}^{\kappa}(A)}:=\|f\|_{L^2(A)}+
\Big(\iint_{\ern\times\ern}|f(x)-f(y)|^2\kappa(x-y)\,\ud y\,\ud x\Big)^{\frac 1 2}.
\end{equation*}
Then the embedding $\mathscr W^\kappa(A)\hookrightarrow L^2(A)$ is compact.
\end{thm}
\begin{thm}\label{compord1}
 Let
$\{s_k\}_{k\in\N}\subset(0,1)$ be a vanishing sequence and let
$\{E_k\}_{k\in\N}\subset\adm$. If
\begin{equation}\label{enbouuno}
\sup_{k\in\N}\Big(\Per^{s_k}_{G}(E_k)-\frac{1}{s_k}\mu^{(s_k)}_{G}(|E_k\cap\Om|)\Big)<+\infty,
\end{equation}
then, up to a subsequence, 
\begin{equation}\label{conveinsi}
\lim_{k\to  +\infty}|E_k\triangle E|=0
\end{equation}
for some set
$E\in\adm$.
\end{thm}
\begin{proof}
Let~$R>1$ be such that~$\Om \subset B_R$.
In view of \eqref{usoora2}, the bound \eqref{enbouuno} remains true  if
{$\mu^{(s_k)}_{G}(|E_k \cap \Om|)$} is replaced by
{$\mu^{(s_k)}_{R,G}(|E_k \cap \Om|)$}. 
{By \eqref{enbouuno}, \eqref{nuovoper1}, \eqref{stimaalto1} and \eqref{stimaalto2}, we have that
\begin{equation}\label{enbouuno+}
+\infty>\sup_{k\in\N}
\Big(\Per^{s_k}_G(E_k)-\frac{1}{s_k}\mu^{(s_k)}_{R,G}(|E_k\cap\Omega|)\Big)
\ge \sup_{k\in\N}\I^{s_k}(E_k\cap\Omega,\Cc E_k\cap\Omega)-C(N,R,\Omega).
\end{equation}
}
{For every $k\in\N$, set $f_k:=\chi_{E_k\cap\Omega}$. Then, setting
$$
\kappa(z):=\frac{\chi_{B_1}(z)}{|z|^N},
$$
 by \eqref{enbouuno+} and Remark \ref{lipsets}, we obtain
\begin{equation*}
\begin{aligned}
&\, \sup_{k\in\N}\frac 1 2\iint_{\R^{N}\times\R^N}|f_{k}(x)-f_{k}(y)|^2\kappa(x-y)\,\ud y\,\ud x
\\
\leq & \,\sup_{k\in\N}\int_{E_{k} \cap \Om}\ud x\int_{\Cc E_{k}\cap\Omega\cap B_1(x)}\frac{\ud y}{|x-y|^N} + \sup_{k\in\N}\int_{E_k \cap \Om}\ud x \int_{\Cc \Om \cap B_1(x)}\frac{\ud y}{|x-y|^N}\\
\le&\,\sup_{k\in\N}\I^{s_k}(E_k\cap\Omega,\Cc E_k\cap\Omega)+\int_{\Omega}\ud x\int_{\Cc \Omega\cap B_{R+1}}\frac{\ud y}{|x-y|^{N}}<+\infty,
\end{aligned}
\end{equation*} 
where in the last 
 inequality we used that~$B_1(x) \subset B_{R+1}$ for any~$x \in \Om \subset B_R$
and that, for every $k\in\N$,
\begin{equation*}
\int_{E_{k} \cap \Om}\ud x\int_{\Cc E_{k}\cap\Omega\cap B_1(x)}\frac{\ud y}{|x-y|^N} \le \int_{E_{k} \cap \Om}\ud x\int_{\Cc E_{k}\cap \Omega}\frac{\ud y}{|x-y|^{N+s_k}}\le \I^{s_k}(E_k\cap\Omega,\Cc E_k\cap\Omega). 
\end{equation*}
Hence, 
since the kernel $\kappa$ is radially symmetric and satisfies~\eqref{iposulkernel}, by Theorem \ref{JaWe} applied with 
$A=\Om$ there exist a (not relabeled) 
subsequence and a function $f \in L^2(\Om)$ such that
$f_{k}\to  f$ strongly in $L^2(\Om)$, so that 
$f=\chi_{F}$ for some $F \in\M(\ern)$ with $F \subset\Om$.  Setting $E=F \cup G$, {we have \eqref{conveinsi}}.
 }
\end{proof}
\subsection{First order $\Gamma$-limit}
We first introduce the candidate $\Gamma$-limit.

To this end, recalling the definition of $\I^s$ in \eqref{intera}, for every $R>0$ and for every $E\in\M(\R^N)$ we set
\begin{eqnarray} 
\label{tilde}
\widetilde{\Per}^0(E)&:=&\,\phantom{+}\I^0(E\cap \Om,\Cc E\cap\Om),\\ 
\label{cappuccio}
\widehat{\Per}^0_{R,G}(E)&:=&\,-|E\cap \Om|\int_{\Om\cap \Cc B_R}\frac{\ud y}{|y|^N}\\ \nonumber
&&\,+\I^{0}(E\cap\Om,\Cc G\cap\Cc \Om\cap B_R)
+\I^0(\Cc E\cap\Om,G\cap B_R)
\\ 
\nonumber
&&\,+\int_{E\cap\Om}\ud x\int_{\Cc G\cap\Cc \Om\cap\Cc B_R}\Big(\frac{1}{|x-y|^N}-\frac{1}{|y|^N}\Big)\,\ud y\\ \nonumber
&&\,+\int_{\Cc E\cap \Om}\ud x\int_{ G\cap \Cc B_R}\Big(\frac{1}{|x-y|^N}-\frac{1}{|y|^N}\Big)\,\ud y.
\end{eqnarray}
Notice that for every $E\in\M(\R^N)$
\begin{equation*}\label{solodentro}
\widetilde{\Per}^0(E)=\widetilde{\Per}^0(E\cap\Om)\quad\textrm{ and }\quad\widehat{\Per}^0_{R,G}(E)=\widehat{\Per}^0_{R,G}(E\cap\Om).
\end{equation*}
Moreover, for every $R>0$ and for every $E\in\M(\R^N)$
\begin{equation}\label{limitato}
|\widehat{\Per}^0_{R,G}(E)|\le M_R(\Om),
\end{equation}
where we have set
\begin{equation}\label{limitato2}
M_R(\Om):=|\Om|\int_{\Om\cap \Cc B_R}\frac{\ud y}{|y|^N}
+{\I^{0}(\Om,\Cc \Om\cap B_R)}
+{\int_{\Om}\ud x\int_{\Cc \Om\cap\Cc B_R}\bigg|\frac{1}{|x-y|^N}-\frac{1}{|y|^N}\bigg|\,\ud y}.
\end{equation}

\begin{rem}\label{continrem}
For every $R>0$ the functional $\widehat{\Per}^0_{R,G}$ is continuous with respect to the strong~$L^1$-convergence in $\Om$; namely, given $\{E_k\}_{k\in\N}\subset\M(\R^N)$ and~$E\in\M(\R^N)$ with
$$
|(E_k\triangle E)\cap\Om|\to0
\qquad\textrm{as }k\to+\infty,
$$
we have
\begin{equation}\label{contpk}
\lim_{k\to +\infty}
\big|\widehat{\Per}^0_{R,G}(E_k)-\widehat{\Per}^0_{R,G}(E)\big|=0.
\end{equation}
Indeed, for every~$E,E'\in\M(\R^N)$, we have
\begin{equation*}\label{contincapp}
\begin{aligned}
\big|\widehat{\Per}^0_{R,G}(E)-\widehat{\Per}^0_{R,G}(E')\big|
\le&\, |(E\triangle E')\cap \Om|\int_{\Om\cap \Cc B_R}\frac{\ud y}{|y|^N}\\
&\,+\I^{0}((E\triangle E')\cap\Om,\Cc G\cap\Cc \Om\cap B_R)\\
&\,+\I^0((E\triangle E')\cap\Om,G\cap B_R)\\
&\,+\int_{(E\triangle E')\cap\Om}\ud x
\int_{\Cc G\cap\Cc \Om\cap\Cc B_R}
\Big|\frac{1}{|x-y|^N}-\frac{1}{|y|^N}\Big|\,\ud y\\
&\,+\int_{(E\triangle E')\cap \Om}\ud x
\int_{G\cap \Cc B_R}
\Big|\frac{1}{|x-y|^N}-\frac{1}{|y|^N}\Big|\,\ud y,
\end{aligned}
\end{equation*}
whence \eqref{contpk} follows  by the absolute continuity of the Lebesgue integral, since all the  {integrands are in $L^1$}
on the corresponding product domains.

Moreover, if $G,G'\in\M(\R^N)$ with~$G,G'\subset\Cc\Om$, then, for every~$E\in\M(\R^N)$,
\begin{equation}\label{contincappG}
\begin{aligned}
\big|\widehat{\Per}^0_{R,G}(E)-\widehat{\Per}^0_{R,G'}(E)\big|
\le&\,\phantom{+} \I^{0}(E\cap\Om,(G\triangle G')\cap\Cc \Om\cap B_R)\\
&\,+\I^0(\Cc E\cap\Om,(G\triangle G')\cap B_R)\\
&\,+\int_{E\cap\Om}\ud x
\int_{(G\triangle G')\cap\Cc \Om\cap\Cc B_R}
\Big|\frac{1}{|x-y|^N}-\frac{1}{|y|^N}\Big|\,\ud y\\
&\,+\int_{\Cc E\cap \Om}\ud x
\int_{(G\triangle G')\cap \Cc B_R}
\Big|\frac{1}{|x-y|^N}-\frac{1}{|y|^N}\Big|\,\ud y.
\end{aligned}
\end{equation}
 Consequently, if $|G_k\triangle G|\to0$, with~$G_k,G\subset\Cc\Om$, then
$$
\lim_{k\to+\infty} \big|\widehat{\Per}^0_{R,G_k}(E)-\widehat{\Per}^0_{R,G}(E)\big|=0  \qquad\textrm{for every }E\in\M(\R^N).
$$
Since the above estimates are uniform, we also obtain that, whenever
$$
|(E_k\triangle E)\cap\Om| \to 0 \qquad\textrm{and}\qquad |G_k\triangle G| \to 0,
$$
one has
$$
\widehat{\Per}^0_{R,G_k}(E_k)\to \widehat{\Per}^0_{R,G}(E).
$$
\end{rem}


The candidate $\Gamma$-limit $\Per^0_{G}$ of the functionals $\Per^{s_k}_{G}-\frac{1}{s_k}\mu^{(s_k)}_{G}$ is defined for every $E\in\M(\R^N)$ as
{ 
\begin{equation}\label{per1}
\Per^0_{G}(E) :=
\begin{cases}
\widetilde{\Per}^0(E)+\widehat{\Per}^0_{1,G}(E) & \text{if}~E \in \adm,\\
+\infty & \text{otherwise}.
\end{cases}
\end{equation}
}
In view of \eqref{limitato}, the functional $\Per^0_G$ is well-defined and takes values in $[-M_1(\Om),+\infty]$.

The next result - proven in Appendix \ref{appendice} - shows that, for $G=\emptyset$, $\Per^0_G$ agrees with  
 the $0$-fractional perimeter $\hat H^0$ defined in \cite{DLNP21} for finite volume sets $E$ as
\begin{equation}\label{defH0}
\hat H^0(E):=\int_{E}\ud x\int_{\Cc E\cap B_1(x)}\frac{\ud y}{|x-y|^N}-\int_{E}\ud x\int_{E\cap\Cc B_1(x)}\frac{\ud y}{|x-y|^N}.
\end{equation}
\begin{prop}\label{comezeroper}
 For every $E\in\mathscr{A}_{\emptyset}$ we have
\begin{equation*}\label{vuoto}
\Per^0_{\emptyset}(E)=\hat H^0(E).
\end{equation*}
\end{prop}
Next, we observe how the functional $\widehat\Per^0_{R,G}$ changes with respect to $R$. To this end, for every $0<R_1<R_2$ we set
$A_{R_1,R_2}:=B_{R_2}\setminus B_{R_1}$.
Moreover, for every $F\in\M(\ern)$, 
we set
\begin{equation}\label{costadd}
\hat\alpha_{R_1,R_2}^{0}(F):=\int_{F\cap A_{R_1,R_2}}\frac{\ud z}{|z|^N},
\end{equation}
and we observe that
\begin{equation}\label{somma}
\hat\alpha_{R_1,R_2}^{0}(F)+\hat\alpha_{R_1,R_2}^{0}(\Cc F)=N\omega_N\log\frac{R_2}{R_1}.
\end{equation} 
\begin{lemma}
Let $0<R_1<R_2$. For every 
$E\in\M(\R^N)$
\begin{equation}\label{usogamma1}
\begin{aligned}
\,\widehat\Per^{0}_{R_1,G}(E) \,
=\,\widehat\Per^{0}_{R_2,G}(E)-\big(\hat\alpha^{0}_{R_1,R_2}(\Cc G)|E\cap\Om|+\hat\alpha^{0}_{R_1,R_2}(G)|\Cc E\cap\Om|\big).
\end{aligned}
\end{equation} 
\end{lemma}
\begin{proof}
It is straightforward to check that
\begin{equation}\label{primog}
\begin{aligned}
&\,-|E\cap \Om|\int_{\Om\cap \Cc B_{R_1}}\frac{\ud y}{|y|^N}\\
=&\,-|E\cap \Om|\int_{\Om\cap \Cc B_{R_2}}\frac{\ud y}{|y|^N}-|E\cap \Om|\int_{\Om \cap A_{R_1,R_2}}\frac{\ud y}{|y|^N},\\ 
&\,\\
&\,\phantom{+}\I^0(E\cap\Om,\Cc G\cap\Cc \Om\cap B_{R_1})\\
=&\,\phantom{+}\I^0(E\cap\Om,\Cc G\cap\Cc \Om\cap B_{R_2})-\I^0(E\cap\Om,\Cc G\cap\Cc \Om\cap A_{R_1,R_2}),\\ 
&\,\\
&\,\phantom{+}\int_{E\cap \Om}\ud x\int_{\Cc G\cap\Cc\Om\cap \Cc B_{R_1}}\Big(\frac{1}{|x-y|^N}-\frac{1}{|y|^N}\Big)\,\ud y\\
=&\,\phantom{+}\int_{E\cap \Om}\ud x\int_{\Cc G\cap\Cc\Om\cap \Cc B_{R_2}}\Big(\frac{1}{|x-y|^N}-\frac{1}{|y|^N}\Big)\,\ud y\\
&\,+\int_{E\cap \Om}\ud x\int_{\Cc G\cap\Cc\Om\cap A_{R_1,R_2}}\Big(\frac{1}{|x-y|^N}-\frac{1}{|y|^N}\Big)\,\ud y;
\end{aligned}
\end{equation}
moreover,  since $\Cc G=\Om\cup(\Cc G\cap\Cc\Om)$ (where the union is disjoint), by \eqref{costadd}, we have
\begin{equation}\label{primog0}
\begin{aligned}
&\,-|E\cap \Om|\int_{\Om \cap A_{R_1,R_2}}\frac{\ud y}{|y|^N}-\I^0(E\cap\Om,\Cc G\cap\Cc \Om\cap A_{R_1,R_2})\\
&\,+\int_{E\cap \Om}\ud x\int_{\Cc G\cap\Cc\Om\cap A_{R_1,R_2}}\Big(\frac{1}{|x-y|^N}-\frac{1}{|y|^N}\Big)\,\ud y\\
=&\,-\hat{\alpha}^0_{R_1,R_2}(\Cc G)|E\cap\Om|.
\end{aligned}
\end{equation}
Therefore, by \eqref{primog} and \eqref{primog0}, we have
 \begin{equation}\label{primogg}
\begin{aligned}
&\,-|E\cap \Om|\int_{\Om\cap \Cc B_{R_1}}\frac{\ud y}{|y|^N}+\I^0(E\cap\Om,\Cc G\cap\Cc \Om\cap B_{R_1})\\
&\,+\int_{E\cap \Om}\ud x\int_{\Cc G\cap\Cc\Om\cap \Cc B_{R_1}}\Big(\frac{1}{|x-y|^N}-\frac{1}{|y|^N}\Big)\,\ud y\\
=&\,-|E\cap \Om|\int_{\Om\cap \Cc B_{R_2}}\frac{\ud y}{|y|^N}+\I^0(E\cap\Om,\Cc G\cap\Cc \Om\cap B_{R_2})\\
&\,+
\int_{E\cap \Om}\ud x\int_{\Cc G\cap\Cc\Om\cap \Cc B_{R_2}}\Big(\frac{1}{|x-y|^N}-\frac{1}{|y|^N}\Big)\,\ud y\\
&\,-\hat{\alpha}^0_{R_1,R_2}(\Cc G)|E\cap\Om|. 
\end{aligned}
\end{equation}
Analogously,
\begin{equation*}
\I^0(\Cc E\cap \Om,G\cap B_{R_1})=\I^0(\Cc E\cap \Om,G\cap B_{R_2})-\I^0(\Cc E\cap \Om,G\cap A_{R_1,R_2})
\end{equation*}
\begin{equation*}
\begin{aligned}
\int_{\Cc E\cap \Om}\ud x\int_{ G\cap \Cc B_{R_1}}\Big(\frac{1}{|x-y|^N}-\frac{1}{|y|^N}\Big)\,\ud y
=&\,\phantom{+}\int_{\Cc E\cap \Om}\ud x\int_{ G\cap \Cc B_{R_2}}\Big(\frac{1}{|x-y|^N}-\frac{1}{|y|^N}\Big)\,\ud y\\
&\,+\int_{\Cc E\cap \Om}\ud x\int_{ G\cap A_{R_1,R_2}}\Big(\frac{1}{|x-y|^N}-\frac{1}{|y|^N}\Big)\,\ud y,
\end{aligned}
\end{equation*}
so that
\begin{equation}\label{secondogg}
\begin{aligned}
&\,\phantom{+}\I^0(\Cc E\cap \Om,G\cap B_{R_1})+\int_{\Cc E\cap \Om}\ud x\int_{ G\cap \Cc B_{R_1}}\Big(\frac{1}{|x-y|^N}-\frac{1}{|y|^N}\Big)\,\ud y\\
=&\,\phantom{+}\I^0(\Cc E\cap \Om,G\cap B_{R_2})+\int_{\Cc E\cap \Om}\ud x\int_{ G\cap \Cc B_{R_2}}\Big(\frac{1}{|x-y|^N}-\frac{1}{|y|^N}\Big)\,\ud y\\
&\,-\hat\alpha^0_{R_1,R_2}(G)|\Cc E\cap\Om|.
\end{aligned}
\end{equation}
By \eqref{cappuccio}, summing \eqref{primogg} and \eqref{secondogg}, we get \eqref{usogamma1}.
\end{proof}
\begin{rem}\label{rem:neiminimi}
For every $E\in\mathscr{A}_{\emptyset}$ we have $\Cc E\in\mathscr{A}_{\Cc\Om}$ and
\begin{equation}\label{neiminimi}
\Per^0_{\emptyset}(E)=\Per^{0}_{\Cc\Om}(\Cc E)+\big(|\Cc E\cap\Om|-|E\cap\Om|\big)\int_{\Om\cap\Cc B_1}\frac{\ud z}{|z|^N}.
\end{equation}
Indeed, let $R>0$ with  $\Om\subset B_R$. Then
\begin{equation}\label{dura}
\begin{aligned}
\widehat\Per^0_{R,\emptyset}(E)
=&\,\I^0(E\cap\Om,\Cc\Om\cap B_R)+\int_{E\cap\Om}\ud x\int_{\Cc B_R}\Big(\frac{1}{|x-y|^N}-\frac{1}{|y|^N}\Big)\ud y\\
=&\,\widehat\Per^0_{R,\Cc\Om}(\Cc E),
\end{aligned}
\end{equation}
where in the last equality we used 
\begin{equation*}
\begin{aligned}
\I^0(E\cap\Om,\Cc\Om\cap B_R)=&\,\I^0(\Cc(\Cc E)\cap\Om,\Cc\Om\cap B_R),\\
\int_{E\cap\Om}\ud x\int_{\Cc B_R}\Big(\frac{1}{|x-y|^N}-\frac{1}{|y|^N}\Big)\ud y
=&\,\int_{\Cc(\Cc E)\cap\Om}\ud x\int_{\Cc\Om\cap\Cc B_R}\Big(\frac{1}{|x-y|^N}-\frac{1}{|y|^N}\Big)\ud y.
\end{aligned}
\end{equation*}
Therefore, enforcing $R>1$, by   \eqref{usogamma1}, \eqref{dura}, and \eqref{somma}, we have
\begin{equation*}
\begin{aligned}
\widehat\Per^0_{1,\emptyset}(E)=&\,\widehat\Per^0_{R,\emptyset}(E)-N\omega_N\log R|E\cap\Om|\\
=&\,\widehat\Per^0_{R,\Cc \Om}(\Cc E)-N\omega_N\log R|E\cap\Om|\\
=&\,\widehat\Per^0_{1,\Cc\Om}(\Cc E)+\hat\alpha^0_{1,R}(\Om)|\Cc E\cap\Om|+\hat\alpha^0_{1,R}(\Cc\Om)|E\cap\Om|-N\omega_N\log R|E\cap\Om|\\
=&\,\widehat\Per^{0}_{1,\Cc\Om}(\Cc E)+\big(|\Cc E\cap\Om|-|E\cap\Om|\big)\int_{\Om\cap\Cc B_1}\frac{\ud z}{|z|^N},
\end{aligned}
\end{equation*}
whence, {using that 
\begin{equation}\label{alcompl}
\widetilde{\Per}^0(E)=\widetilde{\Per}^0(\Cc E)\qquad\textrm{for every }E\in\M(\R^N),
\end{equation}}
\eqref{neiminimi} follows by the very definition of $\Per^0$ in \eqref{per1}. 
\end{rem}
Now we state the first-order $\Gamma$-convergence result.
\begin{thm}\label{thmprimoord}
Let $\{s_k\}_{k\in\N}\subset(0,1)$ be a vanishing sequence. Then, the following $\Gamma$-convergence result holds true.
\begin{itemize}
\item[(i)] ($\Gamma$-liminf inequality) Let $E\in\adm$ and let $\{E_k\}_{k\in\N}\subset\M(\R^N)$ satisfy \eqref{conveinsi}.
Then
\begin{equation}\label{lowerboundprimo}
\Per^0_{G}(E)\le \liminf_{k\to +\infty} \Big(\Per^{s_k}_{G}(E_k)-\frac{1}{s_k}\mu^{(s_k)}_{G}(|E_k\cap\Om|)\Big).
\end{equation}
\item[(ii)] ($\Gamma$-limsup inequality) For every $E\in\adm$ there exists  a sequence $\{E_k\}_{k\in\N}\subset\adm$ 
satisfying \eqref{conveinsi} and
\begin{equation*}
\Per^0_{G}(E)=\lim_{k\to +\infty} \Big(\Per^{s_k}_{G}(E_k)-\frac{1}{s_k}\mu^{(s_k)}_{G}(|E_k\cap\Om|)\Big).
\end{equation*}
\end{itemize}
\end{thm}
We prove (i) and (ii) of Theorem \ref{thmprimoord} in Subsections \ref{subs:lowerbound} and \ref{subs:upperbound}, respectively.

\medskip
\subsection{Lower bound}\label{subs:lowerbound}
In analogy with \eqref{per1}, for every $R>0$ and for every $E\in\M(\R^N)$ we define
{ \begin{equation}\label{perro}
\Per^0_{R,G}(E):=
\begin{cases}
\widetilde\Per^0(E)+\widehat\Per^0_{R,G}(E),
&\textrm{if }E\in\adm,\\
+\infty,
&\textrm{otherwise}.
\end{cases}
\end{equation}}
so that $\Per^0_G\equiv \Per^0_{1,G}$.
Note that, in view of \eqref{usogamma1}, for every $E\in\M(\ern)$ and for every~$R>1$
\begin{equation}\label{usogamma10}
\Per^0_G(E)=\Per^0_{R,G}(E)-\big(\hat\alpha^{0}_{1,R}(\Cc G)|E\cap\Om|+\hat\alpha^{0}_{1,R}(G)|\Cc E\cap\Om|\big).
\end{equation}
In the next remark we rewrite $\Per^0_{R,G}$ in the case that $\Om\subset { B_R}$.
\begin{rem}\label{rmk:perro2}
Let $R>0$ be such that $\Om\subset B_R$. Then, for  every $E\in\adm$
\begin{equation*}
\I^0(E\cap\Om,\Cc E\cap\Om)+\I^0(E\cap\Om,\Cc G\cap\Cc\Om\cap B_R)
=\I^0(E\cap\Om,\Cc E\cap B_R),
\end{equation*}
so that, by \eqref{perro}, recalling \eqref{tilde} and \eqref{cappuccio} we get
\begin{equation}\label{perro2}
\begin{aligned}
\Per^0_{R,G}(E)=&\,\phantom{+}\I^0(E\cap\Om,\Cc E\cap B_R)+\int_{E\cap\Om}\ud x\int_{\Cc G\cap\Cc B_R}\Big(\frac{1}{|x-y|^N}-\frac{1}{|y|^N}\Big)\,\ud y\\
&\,+\I^0(\Cc E\cap\Om,G\cap B_R)+\int_{\Cc E\cap\Om}\ud x\int_{G\cap\Cc B_R}
\Big(\frac{1}{|x-y|^N}-\frac{1}{|y|^N}\Big)\,\ud y.
\end{aligned}
\end{equation}
\end{rem}

We are now in a position to prove the $\Gamma$-liminf inequality in Theorem \ref{thmprimoord}.
\begin{proof}[Proof of Theorem \ref{thmprimoord} (i)] 
We can assume without loss of generality that the bound \eqref{enbouuno} is satisfied, so that $\{E_k\}_{k\in\N}\subset\adm$.
{
By Fatou's Lemma, in view of the definition of $\widetilde\Per^0$ in \eqref{tilde}, we have immediately
\begin{equation}\label{zeropezzouno}
\liminf_{k\to +\infty}\I^{s_k}(E_k\cap\Omega,\Cc E_k\cap\Omega)\ge \widetilde{\Per}^{0}(E).
\end{equation}}
Let moreover { $R>1$} be such that $\Om\subset B_R$.
Then, by \eqref{crucc1},
using {the Dominated Convergence Theorem},  we have
\begin{equation}\label{primopezzouno}
\begin{aligned}
&\,\displaystyle{\lim_{k\to +\infty}}\Big(\I^{s_k}(E_k\cap\Om,\Cc G\cap\Cc\Omega)-\frac{\alpha^{s_k}_{R}(\Cc G)}{s_k}|E_k\cap\Om|\Big)\\
=&\,\I^0(E\cap\Om, \Cc G \cap \Cc \Om \cap B_R)+\int_{E\cap\Om}\ud x\int_{\Cc G\cap\Cc B_R}\Big(\frac{1}{|x-y|^N}-\frac{1}{|y|^N}\Big)\,\ud y,
\end{aligned}
\end{equation}
and, analogously, {by \eqref{crucc2}},
\begin{equation}\label{secondopezzouno}
\begin{aligned} 
&\,{{\lim_{k\to +\infty}}}\Big(\I^{s_k}(\Cc E_k\cap\Om,G\cap\Cc\Om)-\frac{\alpha^{s_k}_{R}(G)}{s_k}|\Cc E_k\cap\Om|\Big)\\
=&\,\I^0(\Cc E\cap\Om,G\cap B_R)+\int_{\Cc E\cap\Om}\ud x\int_{G\cap\Cc B_R}\Big(\frac{1}{|x-y|^N}-\frac{1}{|y|^N}\Big)\,\ud y.
\end{aligned}
\end{equation}

Moreover, by the first equality in \eqref{usoora} and by \eqref{conveinsi}, we have that
\begin{multline}\label{oraora}
\lim_{k\to +\infty}\frac{1}{s_k}(\mu^{(s_k)}_{R,G}(|E_k\cap\Om|)-\mu^{(s_k)}_{G}(|E_k\cap\Om|))\\
=-\lim_{k\to +\infty}\int_{\Cc G\cap A_{1,R}}\frac{\ud z}{|z|^{N+s_k}}|E_k\cap\Om|-\lim_{k\to +\infty}\int_{G\cap A_{1,R}}\frac{\ud z}{|z|^{N+s_k}}|\Cc E_k\cap\Om|\\
=-\big(\hat\alpha_{1,R}^0(\Cc G)|E\cap\Om|+\hat\alpha_{1,R}^0(G)|\Cc E\cap\Om|\Big),
\end{multline}
where $\hat\alpha_{1,R}^0$ is defined in \eqref{costadd}.
Finally,  by \eqref{nuovoper1}, 
 \eqref{oraora}, \eqref{perro2},   
summing {\eqref{zeropezzouno}}, \eqref{primopezzouno} and \eqref{secondopezzouno}, and using \eqref{usogamma10},  we get
\begin{equation*} 
\begin{aligned}
\liminf_{k\to +\infty}\Big(\Per^{s_k}_{G}(E_k)-\frac{1}{s_k}\mu^{(s_k)}_{G}(|E_k\cap\Om|)\Big)
\ge&\,\phantom{-}\liminf_{k\to +\infty}\Big(\Per^{s_k}_{G}(E_k)-\frac{1}{s_k}\mu^{(s_k)}_{R,G}(|E_k\cap\Om|)\Big)\\
&\,+\lim_{k\to +\infty}\frac{1}{s_k}\big(\mu^{(s_k)}_{R,G}(|E_k\cap\Om|)-\mu^{(s_k)}_{G}(|E_k\cap\Om|)\big)\\ 
\ge&\,\phantom{-}\Per^0_{R, G}(E)-\big(\hat\alpha_{1,R}^0(\Cc G)|E\cap\Om|+\hat\alpha_{1,R}^0(G)|\Cc E\cap\Om|\big)\\
=&\,\phantom{-}\Per^0_{G}(E),
\end{aligned}
\end{equation*}
i.e., \eqref{lowerboundprimo}.
\end{proof}

\medskip
\subsection{Upper bound}\label{subs:upperbound}
Here we prove the $\Gamma$-limsup inequality (ii) of Theorem \ref{thmprimoord}. To this end, we first prove a density result for $\Per^0_{G}$. Specifically, we show that
measurable subsets of $\Om$ with finite relative $0$-fractional perimeter $\Per^0_{G}$ can be approximated in measure and in ($\Per^0_G$)-perimeter by
smooth sets compactly contained in $\Om$.
We will prove such a density result following along the lines of \cite[Theorem 3.42]{AFP00} (see also \cite{Lom18,DLNP21}).
We first show that the nonlocal perimeter $\widetilde\Per^0$ is submodular. 
\begin{lemma}\label{nonlocalper}
The functional $\widetilde\Per^0$ in \eqref{tilde} is sub-modular in $\M(\R^N)$, i.~\!e., for every $E^1,E^2\in\M(\R^N)$
\begin{equation}\label{actsubadd}
\widetilde\Per^0(E^1\cup E^2)+\widetilde\Per^0(E^1\cap E^2)\le \widetilde\Per^0(E^1)+\widetilde\Per^0(E^2).
\end{equation}
\end{lemma}
\begin{proof}
For every $E\in\M(\R^N)$,
\begin{equation*}\label{reprefo}
\widetilde{\Per}^0(E)
=
\frac12\iint_{\Om\times\Om}
\frac{|\chi_E(x)-\chi_E(y)|}{|x-y|^N}\,\dx\,\dy .
\end{equation*}
For $a,b,u,v\in\{0,1\}$ one has
$$
|\max\{a,b\}-\max\{u,v\}|
+
|\min\{a,b\}-\min\{u,v\}|
\le |a-u|+|b-v|.
$$
Applying the above inequality  with
$a=\chi_{E^1}(x)$, $b=\chi_{E^2}(x)$, $u=\chi_{E^1}(y)$ and
$v=\chi_{E^2}(y)$, and using
$$
\chi_{E^1\cup E^2}=\max\{\chi_{E^1},\chi_{E^2}\},
\qquad
\chi_{E^1\cap E^2}=\min\{\chi_{E^1},\chi_{E^2}\},
$$
we obtain \eqref{actsubadd} after multiplication by $|x-y|^{-N}$ and
integration over $\Om\times\Om$.
\end{proof}

\begin{rem}\label{tvconv}
\rm{
We define the {\it generalized total variation functional} $\TV_{\widetilde\Per^{0}}:L^1(\Om)\to [0,+\infty]$ as
\begin{equation*}
\TV_{\widetilde\Per^{0}}(u):=\int_{\R}\widetilde\Per^0(\{u>t\}) \, \ud t\qquad\textrm{for every }u\in L^1(\Om).
\end{equation*}
By Lemma \ref{nonlocalper}, in view of \cite[Proposition 3.4]{CGL10}, the functional $\TV_{\widetilde\Per^0}$ is convex.
}
\end{rem}
We have the following density result for $\widetilde{\Per}^0$.
\begin{prop}\label{solounpezz} 
{For every $F\in\M(\R^N)$ with $\mathrm{dist}(F,\Cc\Om)>0$, and {$\widetilde{\Per}^0(F)<+\infty$}, there exists a sequence $\{F_k\}_{k\in\N}$ of smooth sets  such that}
\begin{equation}\label{vincoli}
 \mathrm{dist}(F_k,\Cc\Om)>0,\quad |F_k|=|F|\qquad\textrm{for every }k\in\N,
\end{equation}
and
\begin{equation}\label{conve}
|F_k\triangle F|\to 0\quad \textrm{and}\quad \widetilde{\Per}^0(F_k)\to \widetilde{\Per}^0(F)\qquad\textrm{as }k\to +\infty.
\end{equation}
\end{prop}
\begin{proof} 
If $|F|=0$, we take~$F_k:=\emptyset$ for every $k\in\N$ and there is nothing
left to prove. Hence, we assume that $|F|>0$. 
Set $\delta:=\mathrm{dist}(F,\Cc\Om)>0$.
By the same lower semicontinuity, convexity and coarea arguments used in
\cite[Proposition~7.2]{DLNP21}, applied here to the convex functional~$\TV_{\widetilde\Per^0}$ introduced in Remark~\ref{tvconv}, there exists a sequence~$\{\widetilde F_k\}_{k\in\N}\subset\M(\R^N)$ of smooth sets such that
\begin{equation*}\label{dispo}
\mathrm{dist}(\widetilde F_k,\Cc\Om)>\frac{\delta}{2}\qquad\textrm{for }k\textrm{ large enough},
\end{equation*}
and
$$
|\widetilde F_k\triangle F|\to 0, \qquad \widetilde{\Per}^0(\widetilde F_k)\to \widetilde{\Per}^0(F)
\qquad\textrm{as }k\to+\infty .
$$
Up to discarding finitely many
indices, $|\widetilde F_k|>0$. We set
$$
\lambda_k:=\left(\frac{|F|}{|\widetilde F_k|}\right)^{1/N},
\qquad F_k:=\lambda_k\widetilde F_k .
$$
Then $|F_k|=|F|$ and, since $|\widetilde F_k|\to |F|$, we have~$\lambda_k\to1$.

We now check that the homothetic correction does not spoil neither the constraints nor
the energy convergence. 
Since~$\lambda_k\to1$, for $k$ large enough
$$
|\lambda_k x-x|\le
\frac{\delta}{4} \qquad\textrm{for every }x\in\widetilde F_k.
$$
Therefore
$$
\mathrm{dist}(\lambda_k x,\Cc\Om) \ge \mathrm{dist}(x,\Cc\Om)-|\lambda_k x-x| >\frac{\delta}{4} \qquad\textrm{for every }x\in\widetilde F_k,
$$
and hence, for $k$ large enough
\begin{equation}\label{anchequestamagg}
\mathrm{dist}(F_k,\Cc\Om)>\frac{\delta}{4}, \qquad{\mathrm{dist}(F_k,\Cc(\lambda_k\Om))>\frac{\delta}{4}}.
\end{equation}

Up to relabelling the sequence, this gives the first
condition in~\eqref{vincoli}; the volume constraint has already been built into
the definition of $\lambda_k$.

Moreover, since $\lambda_k \to 1$ as $k\to +\infty$ we have
$$
|F_k\triangle F| \le |\lambda_k\widetilde F_k\triangle \lambda_k F|   +|\lambda_k F\triangle F|
= \lambda_k^N|\widetilde F_k\triangle F|  +|\lambda_k F\triangle F| \to 0,
$$
 namely, the first part of \eqref{conve}.
It remains only to prove the convergence of the internal energy $\widetilde{\Per}^0$. 

Now, by the change of variables $(x,y)\mapsto(\lambda_k x,\lambda_k y)$,
$$
\lambda_k^N\widetilde{\Per}^0(\widetilde F_k)
= \I^0(F_k,\lambda_k\Om\setminus F_k).
$$
As a consequence, using \eqref{anchequestamagg}, for $k$ large enough, we have
$$
\begin{aligned}
\left| \widetilde{\Per}^0(F_k)-\lambda_k^N\widetilde{\Per}^0(\widetilde F_k) \right|
&= \left| \I^0(F_k,\Om\setminus F_k) - \I^0(F_k,\lambda_k\Om\setminus F_k) \right|
\\
&\le \int_{F_k}\ud x \int_{\Om\triangle\lambda_k\Om} \frac{\ud y}{|x-y|^N}\le |\Omega||\Omega\triangle\lambda_k\Omega|\frac{4^N}{\delta^N},
\end{aligned}
$$
which, using again that $\lambda_k\to 1$, yields the second part of \eqref{conve}.
\end{proof}

In Proposition \ref{density} below, we prove that the density result  {above}
extends to sets that are not well-contained in $\Om$.
\begin{prop}\label{density} 
For every $F\in\M(\R^N)$ with $F\subset\Om$, $|\Cc F\cap \Om|>0$ and $\widetilde\Per^0(F)<+\infty$,  there exists a sequence of smooth sets $\{F_k\}_{k\in\N}\subset\M(\R^N)$  such that
 \begin{equation}\label{vincoligen}
 \mathrm{dist}(F_k,\Cc\Om)>0,\quad |F_k|=|F|\qquad\textrm{for every }k\in\N,
 \end{equation}
 and
 \begin{equation}\label{convegen}
 |F_k\triangle F|\to 0\quad \textrm{and}\quad \widetilde{\Per}^0(F_k)\to \widetilde{\Per}^0(F)\qquad\textrm{as }k\to +\infty.
\end{equation}
\end{prop}
\begin{proof}

We claim that, for every $\ep>0$ sufficiently small, there exists a measurable set $F_\ep$ such that
\begin{itemize}
\item[(a)] $\mathrm{dist}(F_\ep,\Cc \Om)\ge\ep$,\vspace{0.8mm}
\item[(b)] $|F|=|F_\ep|$,\vspace{0.5mm}
\item[(c)] $|F\triangle F_\ep|\le \eta_1(\ep)$,\vspace{0.8mm}
\item[(d)] $\displaystyle \big|\widetilde{\Per}^0(F)-\widetilde{\Per}^0(F_\ep)\big|\le{\eta_2}(\ep)$,\vspace{0.8mm}
\end{itemize}
where $\eta_1$ and $\eta_2$ are continuous functions vanishing at $0$.

\medskip
Assuming  that the claim holds true, the statement is a consequence of the following diagonal argument.
Let $\{\ep_k\}_{k\in\N}\subset(0,+\infty)$ be a vanishing sequence and, for every $k\in\N$, let $F_{\ep_k}$ be the set provided  by the claim.  
In view of property (a), Proposition~\ref{solounpezz} yields, for every~$k\in\N$, a sequence of smooth sets~$\{F_{k,h}\}_{h\in\N}$ satisfying
\begin{equation*}
 \mathrm{dist}(F_{k,h},\Cc\Om)>0,
\quad |F_{k,h}|=|F_{\ep_k}|=|F|\qquad\textrm{for every }h\in\N
\end{equation*} 
and
\begin{equation*} 
|F_{k,h}\triangle F_{\ep_k}|\to0
\qquad\text{and}\qquad
\widetilde{\Per}^0(F_{k,h})\to \widetilde{\Per}^0(F_{\ep_k})
\qquad\text{as }h\to+\infty.
\end{equation*}
Hence we may choose~$h(k)\in\N$ so large that
\begin{equation*}
|F_{k,h(k)}\triangle F_{\ep_k}|\le \ep_k \qquad\text{and}\qquad
\big|\widetilde{\Per}^0(F_{k,h(k)})-\widetilde{\Per}^0(F_{\ep_k})\big|\le \ep_k;
\end{equation*}
as a consequence, the sequence $\{F_{k,h(k)}\}_{k\in\N}$ satisfies \eqref{vincoligen} and \eqref{convegen}.

\medskip
It remains to prove the claim. If $\mathrm{dist}(F,\Cc\Om)>0$, then it is enough to take $F_\ep\equiv F$ for every~$\ep>0$ small enough. Therefore, we may assume from now on that
\begin{equation}\label{distanza}
\mathrm{dist}(F,\Cc\Om)=0.
\end{equation}
Furthermore we assume that $F$ coincides with its Lebesgue representative.

Let $\ep>0$. Since $\Om$ has Lipschitz continuous boundary, there exists a finite family of points~$\{\xi^j_\ep\}_{j=1}^{J_\ep}\subset \partial \Om$
such that
\begin{equation}\label{stimaJ}
J_\ep\le C_\Om \ep^{1-N},
\end{equation}
and
$$
\partial\Om\subset \bigcup_{j=1}^{J_\ep}B_{\ep}(\xi^j_\ep).
$$
Setting $S_\ep:=\bigcup_{j=1}^{J_\ep}B_{2\ep}(\xi^j_\ep)$, we have
that $\{x\in\Om:\dist(x,\Cc\Om) < \ep\}\subset S_\ep$.
Indeed, if $x\in\Om$ and~$\dist(x,\Cc\Om)<\ep$, then there exists~$\xi\in\partial\Om$ such that~$|x-\xi|<\ep$. Choosing~$j\in\{1,\ldots,J_\ep\}$ with~$|\xi-\xi^j_\ep|<\ep$,
we get~$|x-\xi^j_\ep|<2\ep$, and hence~$x\in S_\ep$.

Let $\widetilde F_\ep:=F\cap\Cc S_\ep$. By the  inclusion above,
$$
\dist(\widetilde F_\ep,\Cc\Om) \ge \ep,
$$
and therefore $\widetilde F_\ep$ satisfies condition~{\rm(a)}.

Furthermore, by \eqref{distanza} and~\eqref{stimaJ},
\begin{equation}\label{omegak}
0<\tilde\eta_\ep:=|F\setminus \widetilde F_\ep| \le |S_\ep|
 \leq  C_\Om\ep^{1-N}\omega_N(2\ep)^N =: C({\Om,N})\ep.
\end{equation}

Since $|\Om\setminus F|>0$,  
 we may
choose a point $\bar\xi\in\Om\setminus F$. Fix~$\bar r>0$ such that
\begin{equation}\label{def_barr}
B_{2\bar r}(\bar\xi)\subset\Om,
\end{equation}
and let
$$
\bar\delta:=|B_{\bar r}(\bar\xi)\setminus F|.
$$
For every $r\in[0,\bar r]$ let
$g(r):=|B_r(\bar\xi)\setminus F|$.
Recalling that $F$ coincides with its Lebesgue representative, we have that $g(r)>0$ for every $r>0$. 

Moreover,
by definition, $g$ is continuous and non-decreasing,~$g(\bar r)=\bar\delta$ and~$g(0)=0$; therefore, taking~$\ep<\bar\delta/{C(\Om,N)}
$,  {by \eqref{omegak}, $\tilde\eta_\ep<\bar\delta$, so that} there exists~$r_\ep\in(0,\bar r)$ such that~$g(r_\ep)=\tilde\eta_\ep$.

By construction,
\begin{equation}\label{rag0}
r_\ep\to 0\qquad\textrm{as }\ep\to 0.
\end{equation}
Moreover, by \eqref{def_barr}, every point of $B_{\bar r}(\bar\xi)$ has distance at least~$\bar r$ from~$\partial\Om$. Hence, for every~$\ep<\bar r/2$,
\begin{equation}\label{oranum}
B_r(\bar\xi)\cap S_\ep = \emptyset \qquad\textrm{for every }0<r<\bar r.
\end{equation}
In particular,
\begin{equation}\label{errekappa}
|B_{r_\ep}(\bar\xi)\setminus \widetilde F_\ep| = |B_{r_\ep}(\bar\xi)\setminus F| = \tilde\eta_\ep.
\end{equation}
We now define
\begin{equation*}\label{defEk}
F_\ep:=\widetilde F_\ep\cup B_{r_\ep}(\bar\xi).
\end{equation*}
By construction, using \eqref{oranum}, the set~$F_\ep$ satisfies property~{\rm(a)}.

Furthermore, by \eqref{errekappa} and \eqref{omegak},
\begin{equation*}
|F_\ep|=|\widetilde F_\ep|+|B_{r_\ep}(\bar\xi)\setminus \widetilde F_\ep|
=|\widetilde F_\ep|+\tilde\eta_\ep=|F|,
\end{equation*}
i.e., property (b); moreover, again by \eqref{errekappa} and \eqref{omegak}, we have
\begin{equation*}
|F\triangle F_\ep|
\le |F\setminus \widetilde F_\ep|+|B_{r_\ep}(\bar\xi)\setminus \widetilde F_\ep|
=2\tilde\eta_\ep
\le 2C(\Om,N)
\ep,
\end{equation*}
whence property (c) follows taking~$\eta_1(\ep)=2C(\Om,N)
\ep$.

It remains to prove property~(d).  
By \eqref{actsubadd} and \eqref{alcompl}, we have
\begin{equation}\label{differen}
\begin{aligned}
\widetilde{\Per}^0(F_\ep)\le&\,\widetilde{\Per}^0(\widetilde F_\ep)+\widetilde{\Per}^0(B_{r_\ep}(\bar\xi))\le \widetilde{\Per}^0(F)+\widetilde{\Per}^0(S_\ep)+\widetilde{\Per}^0(B_{r_\ep}(\bar\xi)).
\end{aligned}
\end{equation}
Note that
\begin{equation}\label{duecontr}
\begin{aligned}
\widetilde\Per^0(B_{r_\ep}(\bar\xi))\le&\,\phantom{+}\int_{B_{r_\ep}(\bar\xi)}\ud x\int_{B_{\mathrm{diam}(\Om)}(x)\setminus B_{r_\ep}(\bar\xi)}\frac{\ud y}{|x-y|^N}\\
\le&\,\phantom{+}\int_{B_{r_\ep}(\bar\xi)}\ud x\int_{B_{2r_\ep}(\bar\xi)\setminus B_{r_\ep}(\bar\xi)}\frac{\ud y}{|x-y|^N}\\
&\,+\int_{B_{r_\ep}(\bar\xi)}\ud x\int_{B_{\mathrm{diam}(\Om)}(x)\setminus B_{r_\ep}(x)}\frac{\ud y}{|x-y|^N}\\
=&\,\int_{B_{r_\ep}(\bar\xi)}\ud x\int_{B_{2r_\ep}(\bar\xi)\setminus B_{r_\ep}(\bar\xi)}\frac{\ud y}{|x-y|^N}+N\omega_N^2r_\ep^N \log\frac{\mathrm{diam(\Om)}}{r_\ep}.
\end{aligned}
\end{equation}
Moreover, by the change of variables $x=\bar\xi+r_\ep x'$, $y=\bar\xi+r_\ep y'$, we have
\begin{equation*}\label{primocontr}
\begin{aligned}
&\,\int_{B_{r_\ep}(\bar\xi)}\ud x\int_{B_{2r_\ep}(\bar\xi)\setminus B_{r_\ep}(\bar\xi)}\frac{\ud y}{|x-y|^N}=r_\ep^N\int_{B_{1}}\ud x'\int_{B_{2}\setminus B_{1}}\frac{\ud y'}{|x'-y'|^N}
=Cr_\ep^N, 
\end{aligned}
\end{equation*}
which, together with \eqref{duecontr} and \eqref{rag0}, yields
\begin{equation}\label{solapalla}
\widetilde\Per^0(B_{r_\ep}(\bar\xi))\le Cr_\ep^N\log\frac{1}{r_\ep}.
\end{equation} 
Similarly, for every $j=1,\ldots,J_\ep$, we have
\begin{equation*}
\widetilde{\Per}^0(B_{2\ep}(\xi^j_\ep))\le C\ep^N\log\frac{1}{\ep},
\end{equation*}
so that,  by \eqref{actsubadd} and \eqref{stimaJ}, we obtain
\begin{equation}\label{persubadd}
\widetilde{\Per}^0(S_\ep)
\, \le\, \sum_{j=1}^{J_\ep} \widetilde{\Per}^0(B_{2\ep}(\xi^j_\ep))
\, \le \, {C}\ep\log\frac 1 \ep.
\end{equation}
By \eqref{differen}, \eqref{solapalla} and \eqref{persubadd}, we have
\begin{equation}\label{unverso}
\limsup_{\ep\to 0}\widetilde{\Per}^0(F_\ep)\le \widetilde{\Per}^0(F).\end{equation}
Now, by lower semicontinuity and by property (c),
\begin{equation}\label{altroverso}
\widetilde{\Per}^0(F)\le \liminf_{\ep\to 0}\widetilde{\Per}^0(F_\ep).
\end{equation}
By combining \eqref{unverso} and \eqref{altroverso}, we get property (d), concluding the proof of the claim and of the whole proposition.
\end{proof}

By Proposition \ref{density}, using Remark \ref{continrem}, we have the following density result for the functional $\Per^0_{G}$.
\begin{corollary}\label{densitycor}
For every $E\in\adm$ with $|\Cc E\cap \Om|>0$,  there exists a sequence of sets $\{E_k\}_{k\in\N}\subset\adm$ with $E_k\cap\Om$ smooth such that
\begin{equation*}\label{vincoligencor}
\mathrm{dist}(E_k\cap\Om,\Cc\Om)>0,\quad |E_k\cap\Om|=|E\cap\Om|\qquad\textrm{for every }k\in\N,
\end{equation*}
and
\begin{equation*}\label{convegencor}
|E_k\triangle E|\to 0\quad \textrm{and}\quad {\Per}^0_{G}(E_k)\to {\Per}^0_{G}(E)\qquad\textrm{as }k\to +\infty.
\end{equation*}
\end{corollary}
\medskip

\medskip

 Before proving the $\Gamma$-limsup inequality Theorem \ref{thmprimoord} (ii), we show the pointwise convergence $\Per^{s}_{G}(\cdot)-\frac 1 s\mu^{(s)}_{G}(|\cdot\cap\Om|)$ to $\Per^{0}_{G}(\cdot)$ as $s\to 0^+$ for Lipschitz sets.
\begin{prop}\label{pointwiseprimo}
Let~$E\in\adm$ be such that $E\cap\Om$ has smooth boundary and $\mathrm{dist}(E\cap\Om, \Cc\Om)>0$.
Then
\begin{equation}\label{puntuale1}
\Per^{0}_{G}(E)=\lim_{s\to 0^+}\Big(\Per^{s}_{G}(E)-\frac 1 s\mu^{(s)}_{G}(|E\cap\Om|)\Big).
\end{equation}
Moreover,
\begin{equation}\label{puntuale2}
\Per^{0}_{G}(\Om\cup G)=\lim_{s\to 0^+}\Big(\Per^{s}_{G}(\Om\cup G)-\frac 1 s{ \mu^{(s)}_{G}(|(\Om\cup G)\cap\Om|)}\Big).
\vspace{3mm}
\end{equation}
\end{prop}
\begin{proof}
We first prove \eqref{puntuale1}. {Since $E\cap\Om$ is smooth and $\mathrm{dist}(E\cap\Om,\Cc\Om)>0$, then the $s$-fractional perimeter of $E\cap\Om$ is bounded (see \cite[formula~(1.4)]{DFPV13});
moreover, it is easy to see that, by Dominated Convergence,
\begin{equation}\label{limite-As}
\lim_{s\to 0^+}\I^{s}(E\cap \Om,\Cc E\cap \Omega) = \widetilde{\Per}^0(E).
\end{equation}
Let now { $R>1$} be such that $\Omega\subset B_R$.
Using Dominated Convergence, by \eqref{crucc1} and \eqref{crucc2}, we get
\begin{equation}\label{primopezzounoub}
\begin{aligned}
&\,\lim_{s\to 0^+}\left(\I^{s}(E\cap\Omega, \Cc G\cap \Cc\Omega)-\frac{1}{s}\alpha^{s}_R(\Cc G)|E\cap\Omega|\right)\\
=&\,{\I^0(E\cap\Omega,\Cc G\cap\Cc\Omega\cap B_R)}+\int_{E\cap\Omega}\ud x\int_{ \Cc G\cap \Cc B_R}\bigg(\frac{1}{|x-y|^{N}}-\frac{1}{|y|^{N}}\bigg)\,\ud y,
\end{aligned}
\end{equation} 
and
\begin{equation}\label{secondopezzounoub}
\begin{aligned}
&\,\lim_{s\to 0^+}\left(\I^s(\Cc E\cap\Omega, G)-\frac{1}{s}\alpha^{s}_R(G)|\Cc E\cap\Omega|\right)\\
=&\,{\I^0(\Cc E\cap\Omega, G\cap B_R)}+\int_{\Cc E\cap\Omega}\ud x\int_{  G\cap \Cc B_R}\bigg(\frac{1}{|x-y|^{N}}-\frac{1}{|y|^{N}}\bigg)\,\ud y
.
\end{aligned}
\end{equation}
}

We now record explicitly the change of renormalization radius for the fixed set
$E$. 
{By \eqref{oraora} applied with $E_k\equiv E$, we get
\begin{equation}\label{oraora-fixed}
\lim_{s\to0^+}
\frac{1}{s} \Big( \mu^{(s)}_{R,G}(|E\cap\Om|) - \mu^{(s)}_{G}(|E\cap\Om|) \Big) 
 = \, -\hat\alpha^0_{1,R}(\Cc G)|E\cap\Om| -\hat\alpha^0_{1,R}(G)|\Cc E\cap\Om| .
\end{equation}
Finally, by \eqref{nuovoper1}, 
summing \eqref{limite-As}, \eqref{primopezzounoub} and~\eqref{secondopezzounoub}, and using
\eqref{oraora-fixed} and~\eqref{usogamma10}, we get
$$
\begin{aligned}
\lim_{s\to0^+} \Big( \Per^{s}_{G}(E) - \frac{1}{s}\mu^{(s)}_{G}(|E\cap\Om|) \Big)  =&\,\phantom{+} \lim_{s\to0^+} \Big(\Per^{s}_{G}(E) - \frac{1}{s}\mu^{(s)}_{R,G}(|E\cap\Om|) \Big) \\
 &\,+ \lim_{s\to0^+} \frac{1}{s} \Big( \mu^{(s)}_{R,G}(|E\cap\Om|) - \mu^{(s)}_{G}(|E\cap\Om|) \Big) \\
=&\,\phantom{+} \Per^0_{R,G}(E) -\big( \hat\alpha^0_{1,R}(\Cc G)|E\cap\Om| + \hat\alpha^0_{1,R}(G)|\Cc E\cap\Om| \big)\\ 
=&\,\phantom{+}\Per^0_G(E).
\end{aligned}
$$
This proves \eqref{puntuale1}.
}

The proof of \eqref{puntuale2} is fully analogous, once noticed that
\begin{equation*} 
\Per^s_{G}(\Om\cup G)-\frac 1 s\mu^{(s)}_{R,G}(|\Om|)={\I^s(\Omega,\Cc G\cap\Cc\Om\cap B_R)}
+\int_{\Om}\ud x\int_{\Cc G\cap\Cc B_R}
\left(\frac{1}{|x-y|^{N+s}}-\frac{1}{|y|^{N+s}}\right)\ud y.
\end{equation*}
\end{proof}

We are now in the position to prove the $\Gamma$-limsup inequality in Theorem \ref{thmprimoord}.

\begin{proof}[Proof of Theorem \ref{thmprimoord}(ii)]
If $E\equiv\Om\cup G$, then, by Proposition \ref{pointwiseprimo}, a recovery sequence is provided by the constant sequence $E_k\equiv E$. 

Assume now that $|\Cc E\cap\Om|>0$. Using the density result in Corollary \ref{densitycor}, {we can assume that $E$ satisfies  $\mathrm{dist}(E\cap\Om,\Cc\Om)>0$}
 and that $E\cap\Om$ is smooth; hence, the claim follows again by Proposition \ref{pointwiseprimo} using a standard diagonal argument in $\Gamma$-convergence.
\end{proof}
This concludes the proof of the whole $\Gamma$-convergence result Theorem \ref{thmprimoord}.
\begin{rem}
By arguing verbatim as in the proof of Theorem \ref{thmprimoord} we have that
for every $R>0$ the functionals $\Per^{s}_{G}(\cdot)-\frac{1}{s}\mu^{(s)}_{R,G}(|\cdot\cap\Om|)$ $\Gamma$-converge (as $s\to 0^+$) to the functional $\Per^0_{R,G}$ defined in \eqref{perro}.
\end{rem}

In some respects, it would be more natural to renormalize the perimeters $\Per^s_G$ by subtracting the zero-order $\Gamma$-limit $\mu_G(|\cdot \cap\Omega|)$, scaled by $1/s$. {This can be done assuming that the  scalar first-order tail correction $\lambda_G$ in \eqref{lamG} exists and lies in $\R$.}

Such a constant
coefficient produces, in the first-order variational limit, the volume bias 
$$
\lambda_G\bigl(|\Cc E\cap\Om|-|E\cap\Om|\bigr),
$$
{as shown by the following result}.

\begin{thm}\label{teolam}
Assume that the limit $\lambda_G$ in \eqref{lamG} exists and is finite. 
 Let $\{s_k\}_{k\in\N}\subset(0,1)$ be a vanishing sequence.

Then, the following $\Gamma$-convergence result holds true.
\begin{itemize}
\item[(i)] ($\Gamma$-liminf inequality) Let $E\in\adm$ and let $\{E_k\}_{k\in\N}\subset\M(\R^N)$ satisfy \eqref{conveinsi}.
Then
\begin{equation*}\label{lowerboundprimolam}
\Per^0_{G}(E)+\lambda_G(|\Cc E\cap\Om|-|E\cap\Om|)
\le \liminf_{k\to +\infty} \Big(\Per^{s_k}_{G}(E_k)-\frac{1}{s_k}\mu_{G}(|E_k\cap\Om|)\Big).
\end{equation*}
\item[(ii)] ($\Gamma$-limsup inequality) For every $E\in\adm$ there exists  a sequence $\{E_k\}_{k\in\N}\subset\adm$ satisfying \eqref{conveinsi} 
such that
\begin{equation*}
\Per^0_{G}(E)+\lambda_G(|\Cc E\cap\Om|-|E\cap\Om|)
=\lim_{k\to +\infty} \Big(\Per^{s_k}_{G}(E_k)-\frac{1}{s_k}\mu_{G}(|E_k\cap\Om|)\Big).
\end{equation*}
\item[(iii)] (Compactness) Let  $\{E_k\}_{k\in\N}\subset\M(\R^N)$. If
\begin{equation*}\label{enbouuno1-lambda}
\sup_{k\in\N}\Big(\Per^{s_k}_{G}(E_k)-\frac{1}{s_k}\mu_{G}(|E_k\cap\Om|)\Big)<+\infty,
\end{equation*}
then, up to a subsequence, \eqref{conveinsi} holds true for some set 
$E\in\adm$.
\end{itemize}
\end{thm}
\begin{proof}
Let $\{E_k\}_{k\in\N}\subset\M(\R^N)$ and let $E\in\M(\R^N)$ satisfy \eqref{conveinsi}. Then, by \eqref{lamG}, 
together with \eqref{alpharcomp} (applied here with $R=1$), and \eqref{eq:CG-alpha},
\begin{equation}\label{richcomp}
\begin{aligned}
&\,\lim_{k\to +\infty}\frac{1}{s_k}\big(\mu_{G}^{(s_k)}(|E_k\cap\Om|)-\mu_{G}(|E_k\cap\Om|)\big)\\
=&\,\lim_{k\to +\infty}\frac{\alpha^{s_k}_1(G)-\alpha(G)}{s_k}\Big(|\Cc E_k\cap\Om|-|E_k\cap\Om|\Big)\\
=&\,\lambda_G(|\Cc E\cap\Om|-|E\cap\Om|);
\end{aligned}
\end{equation}
in other words,  
$\Per^{s_k}_{G}(\cdot)-\frac{1}{s_k}\mu_{G}(|\cdot\cap\Om|)$ 
coincides with the functional $\Per^{s_k}_{G}(\cdot)-\frac{1}{s_k}\mu^{(s_k)}_{G}(|\cdot\cap\Om|)$, up to a continuous perturbation, namely, $\frac{1}{s_k}\big(\mu_{G}^{(s_k)}(|\cdot\cap\Om|)-\mu_{G}(|\cdot\cap\Om|)\big)$.
 Hence, the statement is an immediate consequence of Theorems  \ref{thmprimoord} and \ref{compord1}.
\end{proof}


\medskip
\subsection{$\Gamma$-convergence analysis as $s\to 0^+$ under volume constraint}
Here, we develop the $\Gamma$-convergence analysis for the functional $\Per^s_{G}$ under volume constraints in $\Om$. To this end, for every 
 $ m\in[0,|\Om|]$, 
{we set
\begin{equation}\label{minvincset}
\mathscr{A}_{G,m}:=\{E\in\adm\,:\,|E\cap\Omega|=m\}\,,
\end{equation}
and, for every $s\in [0,1)$,}
we define the functional~$\Per^s_{G,m}:\M(\R^N)\to[-M_1(\Om),+\infty]$  {(see \eqref{limitato2})} as
\begin{equation}\label{volconfun}
\Per^s_{G,m}(E)=\left\{
\begin{array}{ll}
\Per^s_{G}(E)&\textrm{if }~{E\in\mathscr{A}_{G,m}}\,,\\
+\infty&\textrm{otherwise}.
\end{array}
\right.
\end{equation}

\begin{thm}\label{teovol0} 
Let $ m\in[0,|\Om|]$. Let $\{s_k\}_{k\in\N}\subset(0,1)$ be a vanishing sequence.
Then, the following $\Gamma$-convergence result holds true.
\begin{itemize}
\item[(i)] ($\Gamma$-liminf inequality) Let {$E\in\mathscr{A}_{G,m}$} and let $\{E_k\}_{k\in\N}\subset\M(\R^N)$ satisfy \eqref{conveinsi}.
Then
\begin{equation*}\label{lowerboundvol}
\Per^0_{G,m}(E)\le \liminf_{k\to +\infty} \Big(\Per^{s_k}_{G,m}(E_k)-\frac{1}{s_k}\mu^{(s_k)}_{G}(m)\Big).
\end{equation*}
\item[(ii)] ($\Gamma$-limsup inequality) For every {$E\in\mathscr{A}_{G,m}$}
there exists  a sequence {$\{E_k\}_{k\in\N}\subset\mathscr{A}_{G,m}$}
 satisfying \eqref{conveinsi}
 such that
\begin{equation*}
\Per^0_{G,m}(E)=\lim_{k\to +\infty} \Big(\Per^{s_k}_{G,m}(E_k)-\frac{1}{s_k}\mu^{(s_k)}_{G}(m)\Big).
\end{equation*}
\item[(iii)] (Compactness) Let $\{E_k\}_{k\in\N}\subset\M(\R^N)$. If
\begin{equation}\label{enbouvol}
\sup_{k\in\N}\Big(\Per^{s_k}_{G,m}(E_k)-\frac{1}{s_k}\mu^{(s_k)}_{G}(m)\Big)<+\infty,
\end{equation}
then, 
up to a subsequence, 
\eqref{conveinsi} holds true for some set  {$E\in\mathscr{A}_{G,m}$.}
\end{itemize}
\end{thm}
\begin{proof}
By \eqref{volconfun}, for every $s\in(0,1)$
\begin{equation*}
\Per^s_{G}\le \Per^s_{G,m},
\end{equation*}
and hence, the $\Gamma$-liminf inequality~(i) plainly follows from Theorem~\ref{thmprimoord}(i). 
The $\Gamma$-limsup inequality~(ii) follows from the recovery sequence construction in Subsection~\ref{subs:upperbound}, which preserves $|E_k\cap\Om|=|E\cap\Om|=m$.

Finally, we prove the compactness statement. Let~$\{E_k\}_{k\in\N}\subset\M(\R^N)$ satisfy \eqref{enbouvol}. 
Then, by the definition of
$\Per^{s_k}_{G,m}$ in \eqref{volconfun}, we have that {$\{E_k\}_{k\in\N}\subset\mathscr{A}_{G,m}$}
 and hence   \eqref{enbouvol} reads as
$$
\sup_{k\in\N}
\left(\Per^{s_k}_{G}(E_k) -\frac{1}{s_k}\mu^{(s_k)}_{G}{(|E_k\cap\Omega|)}\right)<+\infty .
$$

Therefore Theorem~\ref{compord1} applies, and, up to a subsequence,
\eqref{conveinsi} holds true for some set~$E\in\adm$. Finally, since~$\chi_{E_k\cap\Om}\to\chi_{E\cap\Om}$ in $L^1(\Om)$,
we get that {$E\in\mathscr{A}_{G,m}$.}
\end{proof}
\begin{rem}
 One can easily  combine Theorem \ref{teovol0} and Theorem \ref{teolam}, obtaining that for $G$ as in Theorem \ref{compord1} satisfying \eqref{lamG} and for $m\in[0,|\Om|]$ the functionals
$
\Per^{s}_{G,m}(\cdot)-\frac{1}{s}\mu_{G}(|\cdot\cap\Om|)
$
$\Gamma$-converge (as $s\to 0^+$) to the functional $\Per^0_{G,m}(\cdot)+\lambda_G(|\Om|-2m)$.
\end{rem}


\medskip
\subsection{$\Gamma$-convergence analysis as $s\to 0^+$ for varying exterior data}\label{subs:extdata}
We conclude this section by observing that all the $\Gamma$-convergence statements above can be generalized to the case of varying exterior data. For the zero order $\Gamma$-convergence analysis, the stability has been already observed in Remark \ref{rem:stable}. To discuss stability also for the first order $\Gamma$-convergence results,  consider a sequence
$\{G_k\}_{k\in\N}\subset\M(\R^N)$ with $G_k\subset\Cc\Om$
such that

\begin{equation}\label{conveG}
\lim_{k\to+\infty}|G_k\triangle G|=0
\end{equation}
for some $G\in\M(\R^N)$ with $G\subset\Cc\Om$.
 Theorem \ref{teocambG} below is the generalization of Theorem \ref{thmprimoord} to the case of varying exterior data.

\begin{thm}\label{teocambG}
Let $\{G_k\}_{k\in\N}\subset\M(\R^N)$, with $G_k \subset \Cc \Om$, satisfy 
\eqref{conveG} for some   $G\in\M(\R^N)$. 
Let $\{s_k\}_{k\in\N}\subset(0,1)$ be a vanishing sequence.
Then, the following $\Gamma$-convergence result holds true.
\begin{itemize}
\item[(i)] ($\Gamma$-liminf inequality) Let $E\in\adm$ and let $\{E_k\}_{k\in\N}\subset\M(\R^N)$ satisfy \eqref{conveinsi}.
Then
\begin{equation*}\label{lowerboundvollam}
\Per^0_{G}(E)\le \liminf_{k\to +\infty} \Big(\Per^{s_k}_{G_k}(E_k)-\frac{1}{s_k}\mu^{(s_k)}_{G_k}(|E_k\cap\Om|)\Big).
\end{equation*}
\item[(ii)] ($\Gamma$-limsup inequality) For every $E\in\adm$ there exists  a sequence $\{E_k\}_{k\in\N}\subset\M(\R^N)$ with $E_k\in\mathscr{A}_{G_k}$ (for every $k\in\N$) satisfying \eqref{conveinsi}
 such that
\begin{equation*}
\Per^0_{G}(E)=\lim_{k\to +\infty} \Big(\Per^{s_k}_{G_k}(E_k)-\frac{1}{s_k}\mu^{(s_k)}_{G_k}(|E_k\cap\Om|)\Big).
\end{equation*}
\item[(iii)] Let $\{E_k\}_{k\in\N}\subset\M(\R^N)$. If
\begin{equation}\label{bouen}
\sup_{k\in\N}\Big(\Per^{s_k}_{G_k}(E_k)-\frac{1}{s_k}\mu^{(s_k)}_{G_k}(|E_k\cap\Om|)\Big)<+\infty,
\end{equation}
then, up to a subsequence, \eqref{conveinsi} holds true
for some set~$E\in\adm$.
\end{itemize}
\end{thm}
\begin{proof}
We first record the elementary stability estimates needed below. 
{In view of \eqref{pertufin},} for every $m \in [0,|\Om|]$,
\begin{equation}\label{muvarconv}
\frac{1}{s_k} \big|\mu^{(s_k)}_{R,G_k}(m)-\mu^{(s_k)}_{R,G}(m)\big|
\le |\Om|{R^{-N-s_k}}|G_k\triangle G| \to 0.
\end{equation}
Moreover, by \eqref{conveG} and Dominated Convergence on bounded annuli,
\begin{equation}\label{annvarconv}
\lim_{k\to +\infty}\int_{G_k\cap A_{R_1,R_2}}\frac{\ud z}{|z|^{N+s_k}}
= \int_{G\cap A_{R_1,R_2}}\frac{\ud z}{|z|^N}
\end{equation}
for every $0<R_1<R_2$, and the same convergence holds for the complements.

We prove the lower bound (i). Let $E\in\adm$ and let
$\{E_k\}_{k\in\N}\subset\M(\R^N)$ satisfy~\eqref{conveinsi}. 
{We may assume without loss of generality that the bound \eqref{bouen} is satisfied, so that}
$E_k\in\mathscr A_{G_k}$. 

Fix $R>1$ with $\Om\subset B_R$. The proof of Theorem~\ref{thmprimoord}{\rm(i)} applies with $G_k$ in place of $G$; the terms involving the varying datum pass to the limit by \eqref{conveG},
\eqref{muvarconv}, \eqref{annvarconv}, and Remark~\ref{continrem}. Thus
\begin{equation*}
\begin{aligned}
&\,\liminf_{k\to+\infty}
\Big(\Per^{s_k}_{G_k}(E_k) -\frac{1}{s_k}\mu^{(s_k)}_{G_k}(|E_k\cap\Om|)\Big)\\
\ge&\,\Per^0_{R,G}(E)
-\big(\hat\alpha^0_{1,R}(\Cc G)|E\cap\Om|
+\hat\alpha^0_{1,R}(G)|\Cc E\cap\Om|\big)
=\Per^0_G(E),
\end{aligned}
\end{equation*}
where the last identity follows from~\eqref{usogamma10}. This proves the liminf inequality.

For the limsup inequality, the case $E=\Om\cup G$ is recovered by taking
$E_k=\Om\cup G_k$ and using the pointwise expansion together with
\eqref{conveG},  \eqref{muvarconv} and~\eqref{annvarconv}. If $|\Cc E\cap\Om|>0$, 
{denoting by $\{\widetilde E_k\}_{k\in\N}$ the recovery sequence constructed in the proof of Theorem \ref{thmprimoord}(ii), it is enough to take $E_k:=\big(\widetilde E_k\setminus G\big)\cup G_k$ to get a recovery sequence for the set $E$.}

It remains to prove compactness. Let~$\{E_k\}_{k\in\N}$ satisfy~\eqref{bouen}, so that  $E_k\in\mathscr A_{G_k}$. The proof of Theorem~\ref{compord1}, with $G_k$ in place of $G$ 
gives, after extraction,
$$
{\chi_{E_k\cap\Omega}\to\chi_{E\cap\Omega}
\qquad\textrm{in }L^1(\Omega)}
$$
for some $E\in\M(\R^N)$. Since $E_k\cap\Cc\Om=G_k$, in view of \eqref{conveG}, we have $E\cap\Cc\Om=G$ and
$$
|E_k\triangle E|
\le |(E_k\triangle E)\cap\Om|+|G_k\triangle G|\to0.
$$
Thus \eqref{conveinsi} holds and $E\in\adm$.
\end{proof}

The following result is the analog of Theorem \ref{teovol0} in the case of varying exterior data. 
\begin{thm}\label{teoval0var}
Let $\{G_k\}_{k\in\N}\subset\M(\R^N)$, with~$G_k \subset \Cc \Om$, satisfy 
\eqref{conveG} for some $G\in\M(\R^N)$.
Let $m\in[0,|\Om|]$.
Let $\{s_k\}_{k\in\N}\subset(0,1)$ be a vanishing sequence.
Then, the following $\Gamma$-convergence result holds true.
\begin{itemize}
\item[(i)] ($\Gamma$-liminf inequality) Let {$E\in\mathscr{A}_{G,m}$} and let $\{E_k\}_{k\in\N}\subset\M(\R^N)$ satisfy \eqref{conveinsi}.
Then
\begin{equation}\label{lowerboundvolvar}
\Per^0_{G,m}(E)\le \liminf_{k\to +\infty} \Big(\Per^{s_k}_{G_k,m}(E_k)-\frac{1}{s_k}\mu^{(s_k)}_{G_k}(m)\Big).
\end{equation}
\item[(ii)] ($\Gamma$-limsup inequality) For every {$E\in\mathscr{A}_{G,m}$}  there exists  a sequence $\{E_k\}_{k\in\N}\subset\M(\R^N)$ with {$E_k\in\mathscr{A}_{G_k,m}$ for every $k\in\N$}  satisfying \eqref{conveinsi}, such that
\begin{equation*}
\Per^0_{G,m}(E)=\lim_{k\to +\infty} \Big(\Per^{s_k}_{G_k,m}(E_k)-\frac{1}{s_k}\mu^{(s_k)}_{G_k}(m)\Big).
\end{equation*}
\item[(iii)] (Compactness) Let $\{E_k\}_{k\in\N}\subset\M(\R^N)$. If
\begin{equation}\label{bouenvolvar}
\sup_{k\in\N}\Big(\Per^{s_k}_{G_k,m}(E_k)-\frac{1}{s_k}\mu^{(s_k)}_{G_k}(m)\Big)<+\infty,
\end{equation}
then, up to a subsequence, \eqref{conveinsi} holds true
for some set   {$E\in\mathscr{A}_{G,m}$.}
\end{itemize}
\end{thm}

\begin{proof}
We prove the liminf inequality (i). 
We may assume without loss of generality that 
 {\eqref{bouenvolvar} holds true so that $E_k\in\mathscr{A}_{G_k,m}$ for every $k\in\N$}. 
By Theorem~\ref{teocambG}{\rm(i)} 
\begin{multline*}
\Per^0_{G,m}(E) = \Per^0_G(E)
\le \liminf_{k\to+\infty} \left( \Per^{s_k}_{G_k}(E_k)-\frac{1}{s_k}\mu^{(s_k)}_{G_k}(|E_k\cap\Om|) \right)
\\
= \liminf_{k\to+\infty} \left( \Per^{s_k}_{G_k,m}(E_k)-\frac{1}{s_k}\mu^{(s_k)}_{G_k}(m) \right),
\end{multline*}
which gives~\eqref{lowerboundvolvar}.

As for the limsup inequality (ii), 
 {we highlight that}
the recovery sequence constructed in
Theorem \ref{teocambG}{\rm(ii)} modifies only the interior phase and preserves
the volume in~$\Om$ so that it applies also to this case.

Finally, we show compactness. By definition of $\Per^{s_k}_{G_k,m}$, \eqref{bouenvolvar} implies \eqref{bouen}, so that, up to a subsequence \eqref{conveinsi} holds true. In view of \eqref{conveG}, we have that $E\in\adm$.
Furthermore, by \eqref{bouenvolvar}, $|E_k\cap\Omega|=m$ for every $k$, so that, in view of \eqref{conveinsi}, $|E\cap\Omega|=m$.
\end{proof}
\begin{rem} {We stress that, in view of the choice of the renormalization $\frac{1}{s}\mu^{(s)}_{R,G}$ (defined in \eqref{mur}),  in all the results  of this section - with the only exception of Theorem \ref{teolam} - we never assumed the existence of the limit $\alpha(G)$ in \eqref{eq:alphaG}.}
\end{rem}
\begin{rem}\label{cambionuc}
{Note that the $\Gamma$-convergence analysis of the relative $s$-fractional perimeter with prescribed exterior datum extends with minor notational changes to the case of anisotropic interactions.
Specifically, given an origin-symmetric convex body $\conv\subset\R^N$ and denoting by $|\cdot|_\conv$ the norm having $\conv$ as the unit ball, for every $0<s<1$, one can consider interactions 
 $\I^s_\conv$
 governed by the kernels $\frac{1}{|\cdot|^{N+s}_\conv}$; the corresponding anisotropic $s$-fractional perimeter $\mathrm{Per}_\conv^s$ is thus defined as in \eqref{per} with $\I^s$ replaced by $\I^s_\conv$. As in \eqref{alfaiso}, for every $R>0$ and for every $0<s<1$ we set
\begin{equation}\label{anisoalfa}
\alpha^s_{\conv,R}(G):=s\int_{G\cap\Cc B_R}\frac{\ud z}{|z|_\conv^{N+s}}.
\end{equation}
Assuming that the limit $\alpha_\conv(G):=\lim_{s\to 0^+}\alpha^s_{\conv,1}(G)$ exists, by arguing verbatim as in the proof of Theorem \ref{thm:G-lim}, one has that the $\Gamma$-limit of the functionals $s\Per_{\conv,G}^s$ (with respect to the topology in Theorem \ref{thm:G-lim}) as $s\to 0^+$ is 
\begin{equation*}
\mu_{\conv,G}(m)=\alpha_\conv(\Cc G)m+\alpha_\conv(G)(|\Omega|-m).
\end{equation*}
We highlight that the pointwise limit in the case $G=\emptyset$ is provided in \cite[Theorem 6]{Ludper} (see also \cite[Theorem 2]{Ludnor}).
Furthermore, by defining $\mu^{(s)}_{\conv,G}$ as in \eqref{mur}, according with \eqref{anisoalfa} and
by arguing verbatim as in this section we deduce that the first order $\Gamma$-limit, as $s\to 0^+$, of the renormalized functionals $\Per^{s}_{\conv,G}(\cdot)-\frac{1}{s}\mu^{(s)}_{\conv, G}(|\cdot\cap\Omega|)$ is given by the functional  
\begin{equation*}
\Per^0_{\conv,G}(E):=\widetilde\Per^0_{\conv}(E)+\widehat{\Per}^0_{\conv,1,G}(E),
\end{equation*}
where
\begin{eqnarray*} 
\widetilde{\Per}_\conv^0(E)&:=&\,\phantom{+}\I^0_\conv(E\cap \Om,\Cc E\cap\Om),\\ 
\widehat{\Per}^0_{\conv,R,G}(E)&:=&\,-|E\cap \Om|\int_{\Om\cap \Cc B_R}\frac{\ud y}{|y|_\conv^N}+\I^{0}_\conv(E\cap\Om,\Cc G\cap\Cc \Om\cap B_R)
+\I^0_\conv(\Cc E\cap\Om,G\cap B_R)
\\ 
&&\,+\int_{E\cap\Om}\ud x\int_{\Cc G\cap\Cc \Om\cap\Cc B_R}\Big(\frac{1}{|x-y|_\conv^N}-\frac{1}{|y|_\conv^N}\Big)\,\ud y\\ \nonumber
&&\,+\int_{\Cc E\cap \Om}\ud x\int_{ G\cap \Cc B_R}\Big(\frac{1}{|x-y|_\conv^N}-\frac{1}{|y|_\conv^N}\Big)\,\ud y.
\end{eqnarray*}
}
\end{rem}
%
%

\section{$\Gamma$-convergence analysis as $s\to 1^-$}\label{sec:Gamma1}

{In this section,} we develop the $\Gamma$-convergence analysis of $(1-s)\Per^s_G$ as
$s\to1^-$. For the relative fractional perimeter without fixed exterior
datum, the corresponding limit was proved in \cite{ADM11}: up to a
dimensional constant, one recovers the Euclidean perimeter inside the domain.
{ 
In the present setting the exterior datum is fixed, and the limit perimeter is
taken in $\overline\Om$. We use the compactness and lower-bound results of
\cite{ADM11}, together with the pointwise limit for scaled fractional
Gagliardo seminorms in~\cite{Dav02}.
}

The following compactness statement is \cite[Theorem~1]{ADM11}, rewritten in
our notation.
\begin{thm}\label{compADM}
Let $O$ be a bounded, open set with Lipschitz continuous boundary.
Let $\{s_k\}_{k\in\N}\subset (0,1)$ be a sequence converging to $1$.
Let $\{E_k\}_{k\in\N}\subset\M(\R^N)$ be such that
\begin{equation}\label{bouenADM}
\sup_{k\in\N}\,(1-s_k)\,\I^{s_k}(E_k\cap O',\Cc E_k\cap O')<+\infty
\qquad\textrm{for every }O'\subset\subset O.
\end{equation}
Then, up to a subsequence, there exists a set
$E\in\M(\R^N)$ of locally finite perimeter in $O$ such that
$$
\chi_{E_k}\to\chi_E \qquad\textrm{in }L^1_{\rm loc}(O).
$$
\end{thm}

The next lower bound is an immediate consequence of
\cite[Lemmas~7, 9, 13, 14]{ADM11}.
\begin{prop}\label{lbADM}
Let $O$ be a bounded and open set. Let $\{s_k\}_{k\in\N}\subset (0,1)$ be a sequence converging to~$1$.
Let $F\in\M(\R^N)$ be a set of finite perimeter in $O$, and let
$\{E_k\}_{k\in\N}\subset\M(\R^N)$ be such that
$$
\chi_{E_k}\to\chi_F \qquad\textrm{in }L^1(O).
$$
Then
$$
\omega_{N-1}\mathcal{H}^{N-1}(\partial^\ast F\cap O)
\le \liminf_{k\to +\infty} (1-s_k)\,\I^{s_k}(E_k\cap O,\Cc E_k\cap O).
$$
\end{prop}
 Here and below $\partial^*$ denotes the so called reduced boundary (see \cite{AFP00}) 
Finally, the following pointwise limit follows from \cite[Theorem~1]{Dav02}.
\begin{thm}\label{ubDav}
Let $O$ be a bounded, open set with Lipschitz continuous boundary.
Let $\{s_k\}_{k\in\N}\subset (0,1)$ be a sequence converging to $1$.
Let $F\subset O$ be a set of finite perimeter.
Then
$$
\lim_{k\to +\infty}(1-s_k)\I^{s_k}(F,\Cc F\cap O)
=
\omega_{N-1}\mathcal{H}^{N-1}(\partial^\ast F\cap O).
$$
\end{thm}

\begin{proof}
For every $k\in\N$ we set 
$$
\rho_k(\cdot):=\frac{1}{N\omega_N(2\mathrm{diam}(O))^{1-s_k}}\frac{1-s_k}{|\cdot|^{N+s_k-1}}\chi_{B_{2\mathrm{diam}(O)}}(\cdot).
$$
It is immediate to check that the kernels $\rho_k$ satisfy the assumptions (1) and (2) in \cite{Dav02} and that
\begin{equation*}
(1-s_k)\,\I^{s_k}(F,\Cc F\cap O)=\frac 1 2 N\omega_N(2\mathrm{diam}(O))^{1-s_k}\int_{O}\ud x\int_{O}\frac{|\chi_F(x)-\chi_F(y)|}{|x-y|}\rho_k(x-y)\ud y.
\end{equation*}
Hence, by applying \cite[Theorem 1]{Dav02} with $\Omega=O$ and $u=\chi_F$, we get
\begin{equation*}
\lim_{k\to +\infty }(1-s_k)\,\I^{s_k}(F,\Cc F\cap O)=\frac{N\omega_N}{2} K_{1,N}\mathcal{H}^{N-1}(\partial^\ast F\cap O),
\end{equation*}
where the constant $K_{1,N}$ is defined in \cite[formula (4)]{Dav02} as
\begin{equation*}
K_{1,N}:=\frac{1}{N\omega_N}\int_{\mathbb{S}^{N-1}}|w_1|\,\ud\mathcal{H}^{N-1}(w)=2\,\frac{\omega_{N-1}}{N\omega_N}.
\end{equation*}
\end{proof}

We now state and prove the $\Gamma$-convergence result for $(1-s)\Per^s_G$ as $s\to 1^-$.
\begin{thm}\label{s1}
Assume that $G$ is a set of locally finite perimeter 
 in $\R^N$.
Let $\{s_k\}_{k\in\N}\subset (0,1)$ be a sequence converging to $1$.
Then the following statements hold.
\begin{itemize}
\item[(i)] ($\Gamma$-liminf inequality) Let $E\in\adm$ be a set of locally
finite perimeter and let $\{E_k\}_{k\in\N}\subset\M(\R^N)$ satisfy
\eqref{conveinsi}. Then
\begin{equation}\label{lbs1}
\omega_{N-1}\mathcal{H}^{N-1}(\partial^\ast E\cap\overline\Om)
\le \liminf_{k\to +\infty} (1-s_k)\Per^{s_k}_{G}(E_k).
\end{equation}
\item[(ii)] ($\Gamma$-limsup inequality) For every set $E\in\adm$ of locally
finite perimeter there exists a sequence $\{E_k\}_{k\in\N}\subset\M(\R^N)$
with $E_k\in\mathscr{A}_{G}$ for every $k\in\N$ and satisfying
\eqref{conveinsi} such that
\begin{equation}\label{ubs1}
\omega_{N-1}\mathcal{H}^{N-1}(\partial^\ast E\cap\overline\Om)
=
\lim_{k\to +\infty} (1-s_k)\Per^{s_k}_{G}(E_k).
\end{equation}
\item[(iii)] (Compactness) Let $\{E_k\}_{k\in\N}\subset\M(\R^N)$. If
\begin{equation}\label{bouensuno}
\sup_{k\in\N}\,(1-s_k)\Per^{s_k}_{G}(E_k)<+\infty,
\end{equation}
then, up to a subsequence, \eqref{conveinsi} holds true
for some set~$E\in\adm$ of locally finite perimeter.
\end{itemize}
\end{thm}
\begin{proof}
Let $\{\Om_j\}_{j\in\N}$ be a decreasing sequence of bounded open Lipschitz
sets such that $\Om\subset\subset\Om_{j+1} \subset\subset \Om_{j}$ for every $j\in\N$ and
$\overline\Om=\bigcap_{j\in\N}\Om_j$. 
Set
$$
A_j:=\Om_j\setminus\overline\Om .
$$
Then $A_j$ is bounded, open and Lipschitz. For every $E\in\M(\R^N)$,
\begin{equation}\label{usobilin}
\begin{aligned}
{\rm Per}_{s_k}(E,\Om)=&\,\I^{s_k}(E\cap\Om_j,\Cc E\cap\Om_j)
-\I^{s_k}(E\cap A_j,\Cc E\cap A_j)\\
&\,+\I^{s_k}(E\cap\Om,\Cc E\cap\Cc\Om_j)
+\I^{s_k}(\Cc E\cap\Om,E\cap\Cc\Om_j).
\end{aligned}
\end{equation}
Moreover, setting $\delta_j:=\dist(\Om,\Cc\Om_j)$,
\begin{equation}\label{distanti}
\begin{aligned}
&\,(1-s_k)\,\I^{s_k}(E\cap\Om,\Cc E\cap\Cc\Om_j)
+(1-s_k)\,\I^{s_k}(\Cc E\cap\Om,E\cap\Cc\Om_j)\\
\le&\,|\Om|N\omega_N\frac{1-s_k}{s_k}\delta_j^{-s_k}\to0
\qquad\textrm{as }k\to+\infty .
\end{aligned}
\end{equation}
Since $G$ has locally finite perimeter, Theorem~\ref{ubDav}, applied with
$O=A_j$ and $F=G\cap A_j$, gives
\begin{equation*}
\lim_{k\to+\infty}(1-s_k)\,\I^{s_k}(G\cap A_j,\Cc G\cap A_j)
=
\omega_{N-1}\mathcal H^{N-1}(\partial^\ast G\cap A_j);
\end{equation*}
{in particular, for every sequence $\{E_k\}_{k\in\N}\subset\adm$
\begin{equation}\label{fuoriSuno1}
\lim_{k\to+\infty}(1-s_k)\,\I^{s_k}(E_k\cap A_j,\Cc E_k\cap A_j)
=
\omega_{N-1}\mathcal H^{N-1}(\partial^\ast G\cap A_j).
\end{equation}
}

We prove {\rm(i)}. We may assume that the liminf in \eqref{lbs1} is finite
and argue along a (not relabeled) subsequence realizing it for which
$E_k\in\adm$. By
\eqref{usobilin}, \eqref{fuoriSuno1}, and by Proposition~\ref{lbADM} applied with
$O=\Om_j$ and $F=E$,
\begin{equation}\label{perliminfs1}
\begin{aligned}
\liminf_{k\to+\infty}\,(1-s_k)\Per^{s_k}_G(E_k)
\ge&\,\omega_{N-1}\mathcal H^{N-1}(\partial^\ast E\cap\Om_j)
-\omega_{N-1}\mathcal H^{N-1}(\partial^\ast G\cap A_j)\\
\ge&\,\omega_{N-1}\mathcal H^{N-1}(\partial^\ast E\cap\overline\Om)
-\omega_{N-1}\mathcal H^{N-1}(\partial^\ast G\cap A_j).
\end{aligned}
\end{equation}
Since $\{A_j\}_{j\in\N}$ is a decreasing sequence whose intersection is empty and $G$ has locally finite perimeter,
$$
\mathcal H^{N-1}(\partial^\ast G\cap A_j)\to0.
$$
Letting $j\to+\infty$ in \eqref{perliminfs1}, we obtain \eqref{lbs1}.

In order to prove {\rm(ii)}, take the constant sequence $E_k\equiv E$. From
\eqref{usobilin}, \eqref{distanti}, and Theorem~\ref{ubDav} with
$O=\Om_j$ and $F=E\cap\Om_j$,
$$
\limsup_{k\to+\infty}\,(1-s_k)\Per^{s_k}_G(E)
\le
\omega_{N-1}\mathcal H^{N-1}(\partial^\ast E\cap\Om_j),
$$
which, sending $j\to+\infty$ (i.e., taking the inf over $j\in\N$), yields \eqref{ubs1}.

It remains to prove compactness. By \eqref{bouensuno},
$E_k\in\adm$ for every $k\in\N$. Taking $j=1$ in
\eqref{usobilin}, using the nonnegativity of the last two terms and
\eqref{fuoriSuno1}, we obtain
$$
\sup_{k\in\N}\,(1-s_k)\,\I^{s_k}(E_k\cap\Om_1,\Cc E_k\cap\Om_1)<+\infty.
$$
Hence \eqref{bouenADM} holds with~$O=\Om_1$. 
 By Theorem~\ref{compADM},
there exist a (not-relabeled) subsequence and a set $F$ of locally finite perimeter in~$\Om_1$ such that~$\chi_{E_k}\to\chi_F$  in~$L^1_{\rm loc}(\Om_1)$.
Since~$\Om\subset\subset\Om_1$, the convergence holds in $L^1(\Om)$.  Thus \eqref{conveinsi} holds with $E:=(F\cap\Omega)\cup G\in\adm$. By construction, $E$ has locally finite perimeter in $\R^N$.
\end{proof}
%
%
\medskip
\section{Fractional isoperimetric problems and first-order selection}
\label{sec_minimi}
We now apply the $\Gamma$-convergence results proved in the previous sections
to analyze the asymptotic behaviour, as $s\to0^+$, of some minimum problems
with prescribed exterior datum.
We first consider the unconstrained
problem
\begin{equation}\label{prob}
\min\bigl\{\Per^s_G(E):\ E\in\adm\bigr\}.
\end{equation}
In such a case, as shown in Section \ref{sec:gamma0th}, the leading-order limit as $s\to 0^+$ is purely volumetric; hence, unless the exterior datum is
balanced at infinity, minimizers collapse to one of the two trivial phases
inside~$\Om$. In the balanced  case \eqref{balcase} the leading term is constant, and the
first nontrivial selection is produced by the first-order limit.

We then impose a volume constraint. In that  {case} the leading tail term is
constant on each fixed-volume class, for every exterior datum. Hence the
first-order functional is the first term capable of selecting the position
and shape of the phase in~$\Om$. We close the section by comparing this
tail-driven selection with the local regime $s\to1^-$.

{
Throughout this section we will give as understood
the existence of minimizers of the problems
\begin{equation}\label{forserimane}
\min\{\Per^s_G(E)\,:\,E\in\adm\}\,,
\end{equation}
and 
\begin{equation}\label{forserimane2}
\min\{\Per^s_G(E)\,:\,E\in\mathscr{A}_{G,m}\},
\end{equation}
where $0<s<1$, $G\in\M(\R^N)$ with $G\subset\Cc\Omega$, $m\in[0,|\Omega|]$. We stress that, in view of \eqref{nuovoper1},
 the existence of such minimizers is a consequence of the Direct Method in the Calculus of Variations, once applied the compactness statement Theorem \ref{JaWe} to the kernel $\kappa(\cdot):=\frac{1}{|\cdot|^{N+s}}$, which ensures that every minimizing sequence $\{E_n\}_{n\in\N}$ for the problem \eqref{forserimane} (resp. \eqref{forserimane2}) converges (up to subsequences) to some set $E\in\adm$ (resp. $E\in\mathscr{A}_{G,m}$).
}

\medskip
\subsection{Unconstrained minimizers}

The following result is a direct consequence of the zero-order
$\Gamma$-convergence theorem.

\begin{prop}\label{casodivunmezzo}
Let $G\in\M(\ern)$ with $G\subset\Cc\Om$. Assume that the limit
$\alpha(G)$ in~\eqref{eq:alphaG} exists. Let
$\{s_k\}_{k\in\N}\subset(0,1)$ be such that $s_k\to0^+$. For every
$k\in\N$, let $E_k$ be a minimizer of
\begin{equation}\label{probn}
\min\{\Per^{s_k}_{G}(E):\ E\in\adm\}.
\end{equation}
Then the following facts hold:
\begin{itemize}
\item[(i)] if $\alpha(G)<\displaystyle\frac{N\omega_N}{2}$, then $|E_k\cap\Omega|\to0$ {and $s_k\Per^{s_k}_G(E_k)\to \alpha(G)|\Omega|$} as $k\to +\infty$;
\item[(ii)] if $\alpha(G)>\displaystyle\frac{N\omega_N}{2}$, then $
|\Omega\setminus E_k|\to0$ {and $s_k\Per^{s_k}_G(E_k)\to (N\omega_N-\alpha(G))|\Omega|$} as $k\to +\infty$.
\end{itemize}
\end{prop}

\begin{proof}
By Remark~\ref{trivia}, up to a subsequence,
$
\chi_{E_k\cap\Om}\ud x\weakstar f\ud x
$
for some $f\in L^\infty(\Om;[0,1])$. 
By Theorem \ref{thm:G-lim}, $f$ is a minimizer of the functional $\mu_G(\|\cdot\|_{L^1(\Omega)})$ in $L^\infty(\Omega;[0,1])$.

By
\eqref{gammaorzero},
\begin{equation}\label{dinuovo}
\mu_G(\|f\|_{L^1(\Om)})
=
\alpha(G)|\Om|
+
\big(N\omega_N-2\alpha(G)\big)\|f\|_{L^1(\Om)}.
\end{equation}
Thus the unique minimizer is $f\equiv0$ if
$\alpha(G)<N\omega_N/2$, and it is $f\equiv1$ if
$\alpha(G)>N\omega_N/2$. Hence every weak-star limit point is uniquely
determined  so that the whole sequence $\chi_{E_k\cap\Om}$ weakly star converges to $f$. In  case (i) this gives
$|E_k\cap\Om|\to0$ { and $s_k\Per^{s_k}_G(E_k)\to \mu_G(0)=\alpha(G)|\Omega|$},
whereas in case (ii) it gives
$|\Om\setminus E_k|\to0$ {and $s_k\Per^{s_k}_G(E_k)\to \mu_G(|\Omega|)=(N\omega_N-\alpha(G))|\Omega|$}.
The proof is concluded.
\end{proof}

The only case not selected by the leading-order functional is the balanced
one \eqref{balcase}.
Indeed, by \eqref{dinuovo}, in this case the functional $\mu_G$ is constant
on $L^\infty(\Om;[0,1])$. Therefore the zero-order theory does not
distinguish between different configurations inside~$\Om$. The next
result shows that, if the first-order tail correction $\lambda_G$ exists (see \eqref{lamG}), minimizers are
selected by the first-order $\Gamma$-limit.

\begin{prop}\label{casounmezzo}
Let $G\in\M(\ern)$, with $G\subset\Cc\Om$, satisfy
\eqref{balcase}. 
Assume moreover that the limit $\lambda_G$ in \eqref{lamG} exists and is finite.
Let
$\{s_k\}_{k\in\N}\subset(0,1)$ be such that $s_k\to0^+$. For every~$k\in\N$, let~$E_k$ be a minimizer of~\eqref{probn}. Then, up to a
subsequence,
$$
|E_k\triangle E_\infty|\to0 \qquad\textrm{as }k\to+\infty
$$
for some~$ E_\infty \in\adm$. Moreover, $E_\infty$ is a minimizer in~$\adm$
of the functional
\begin{equation*}\label{bilanciato}
E \, \mapsto \,
\Per^0_G(E) + \lambda_G\big(|\Cc E\cap\Om|-|E\cap\Om|\big)\,,
\end{equation*}
{and
$$
\lim_{k\to +\infty} \Big(\Per^{s_k}_{G}(E_k)-\frac{N\omega_N}{2s_k}|\Omega|\Big)=\Per^0_G(E_\infty) + \lambda_G\big(|\Cc E_\infty\cap\Om|-|E_\infty\cap\Om|\big).
$$
}
\end{prop}

\begin{proof}
By \eqref{balcase} and \eqref{costa}, the quantity
$\mu_G(|E\cap\Om|)$ is independent of $E\in\adm$.
Therefore, for every $k\in\N$, the sets $E_k$ are also minimizers of
\begin{equation*}
\min\Big\{\Per^{s_k}_{G}(E)-\frac{1}{s_k}\mu_G(|E\cap\Omega|)\,:\, E\in\adm\Big\}.
\end{equation*}
The minimality of $E_\infty$ {and the convergence of the minimal energies are then consequences} of the $\Gamma$-convergence
result in Theorem~\ref{teolam}.
\end{proof}
\begin{rem}\label{estremi}
Proposition~\ref{casounmezzo}, together with Theorem~\ref{teolam}, can be extended to the case $\lambda_G=\pm\infty$ in~\eqref{lamG}.
Here we briefly describe such an extension without providing the details of the proof. Assume $\lambda_G=+\infty$.
Setting $\displaystyle \lambda^{(s)}_G:=\frac{\alpha_1^s(G)-\alpha(G)}{s}$, one can show that the functionals
$$
\frac{1}{\lambda^{(s_k)}_G}\Big(\Per^{s_k}_{G}(\cdot)-\frac{1}{s_k}\mu_{G}(|\cdot \cap\Omega|)\Big)\,,
$$
$\Gamma$-converge, as  $s_k\to 0^+$, to the functional $\overline{\Per}:L^\infty(\Omega;[0,1])$ defined by
$$
\overline{\Per}(f) := \|1-f\|_{L^1(\Om)}-\|f\|_{L^1(\Om)} = |\Om|-2\int_\Om f\,{\rm d}x .
$$
Here the topology is that in the zero-order $\Gamma$-convergence result Theorem~\ref{thm:G-lim}, i.e., the weak-star convergence of the densities $\chi_{E_k\cap\Omega}$ to $f$.
Thus the one-sided first-order term does not depend on $G$ and selects the maximal interior density; hence,
the unique minimizer of $\overline{\Per}$ is $f=1$ a.e. in $\Om$.
Consequently, every family of minimizers of~\eqref{prob} converges, in the
topology of Theorem~\ref{thm:G-lim}, to $\Om\cup G$ as $s\to0^+$.

If instead $\lambda_G=-\infty$, by scaling now the renormalized functionals by $-\lambda_G^{(s)}>0$, we obtain the opposite limiting
functional $-\overline{\Per}$. In such a case, the one-sided first-order term selects
the minimal interior density, namely~$f=0$ a.e. in~$\Om$, and every family
of minimizers of~\eqref{prob} converges to~$G$ as~$s\to0^+$.
\end{rem}

\subsection{Volume constraints and first-order selection}\label{sez_volume1}
{ 
{We now focus on the asymptotic behavior of the volume constrained minimum problems \eqref{forserimane2}
for every $m\in[0,|\Omega|]$.
We recall that {the set $\mathscr{A}_{G,m}$ is defined in \eqref{minvincset} and that, for every $s\in(0,1)$, the perimeter $\Per^s_{G,m}$ is defined in~\eqref{volconfun}}.
}
The next result is
a consequence of
Theorem~\ref{teovol0}. It records, in the form needed below, the convergence
of the constrained minimum values and of constrained minimizers.

\begin{prop}\label{prop_volume-min}
{Let $G\in\M(\ern)$ with $G\subset\Cc\Om$, let $m\in[0,|\Omega|]$ and let $\{s_k\}_{k\in\N}$ be such that $s_k\to 0^{+}$.}
For every $k\in\N$, let
$E_k$ be a minimizer of
\begin{equation}\label{semin}
\min\Big\{\Per^{s_k}_{G}(E):\ E\in\mathscr A_{G,m}\Big\}.
\end{equation}
Then, up to a subsequence,
$$
|E_k\triangle E_\infty|\to0
\qquad\textrm{as }k\to+\infty
$$
for some minimizer $E_\infty$ of $\Per^0_G$ in $\mathscr A_{G,m}$, and
\begin{equation}\label{limiteminimim}
\lim_{k\to +\infty}\Big(\Per^{s_k}_{G}(E_k)-\frac{1}{s_k}\mu^{(s_k)}_G(m)
\Big)=
\Per^0_G(E_\infty).
\end{equation}
If, in addition, the limit in~\eqref{lamG} exists and is finite, then
\begin{equation}\label{volume-vale}
\lim_{k\to +\infty}\Big(\Per^{s_k}_{G}(E_k)
-
\frac{1}{s_k}\mu_G(m)\Big)=
\lambda_G(|\Om|-2m)
+\Per^0_G(E_\infty).
\end{equation}
In particular, in the balanced case \eqref{balcase}, 
\begin{equation}\label{eq_volume-bilanc}
\lim_{k\to +\infty}\Big(\Per^{s_k}_{G}(E_k)-
\frac{N\omega_N}{2s_k}|\Omega|\Big)=\lambda_G(|\Om|-2m)
+\Per^0_G(E_\infty).
\end{equation}
\end{prop}
\begin{proof}
The minimality of $E_\infty$ {and the convergence \eqref{limiteminimim} of the minimal energies are consequences of the $\Gamma$-convergence
result in} Theorem~\ref{teovol0}.
Finally, if the limit in~\eqref{lamG} exists, {by \eqref{richcomp}}, we get
$$
\lim_{s\to 0^+}\frac{\mu^{(s)}_G(m)-\mu_G(m)}{s}  = \lambda_G(|\Om|-2m),
$$
which gives~\eqref{volume-vale}. Since in the balanced case~$\mu_G(m)=\frac{N\omega_N}{2}|\Om|$, \eqref{eq_volume-bilanc} follows.
\end{proof}
}
\begin{rem}\label{unicominimo}
{Let  $G\in\M(\R^N)$ and $m\in[0,|\Omega|]$. If the problem
\begin{equation*}
\min\{\Per^0_G(E)\,:\,E\in\mathscr{A}_{G,m}\}
\end{equation*}
admits a unique minimizer, say $E_\ast$, then, by Proposition \ref{prop_volume-min}, every sequence $\{E_k\}_{k\in\N}\subset\M(\R^N)$ of minimizers of \eqref{semin} converges to $E_\ast$. In particular, $E_\infty\equiv E_\ast$ in Proposition \ref{prop_volume-min}.  The same comment applies to the unconstrained case treated in Proposition \ref{casounmezzo}.
}
\end{rem}

{In view of Propositions \ref{casounmezzo} and \ref{prop_volume-min}, when studying the asymptotic behavior of constrained and unconstrained minimum problems,} both the leading tail term $\mu_G(m)$ and the scalar correction
$\lambda_G(|\Om|-2m)$ are constant. Shape selection can only come from the
spatial part of the first-order functional~$\Per^0_G$.
We make this mechanism explicit in the next subsection.

\medskip
\subsection{Balanced windows and holes}\label{subseccone}

In this subsection we assume that
\begin{equation}\label{Omega-in-B1}
\overline\Om\subset B_1.
\end{equation}
Let $A\subset\mathbb S^{N-1}$ be a measurable set satisfying
$\mathcal{H}^{N-1}(A)=\frac{N\omega_N}{2}$.
Set
\begin{equation}\label{def_cono_a}
C_A:=\{r\theta:\ r\ge0,\ \theta\in A\},
\qquad
G_A:=C_A\cap\Cc B_1.
\end{equation}
By polar coordinates, for every $s\in(0,1)$,
\begin{equation}\label{alfa1cono}
\alpha^s_1(G_A)=\mathcal{H}^{N-1}(A)=\frac{N\omega_N}{2}.
\end{equation}
Let $P\subset \Cc G_A\cap\Cc B_1$ and $Q\subset G_A$ be  {two sets with finite measure}; we
define
\begin{equation}\label{def_gpq}
G:=(G_A\setminus Q)\cup P.
\end{equation}
For every measurable set $S\subset\Cc B_1$, we define $\Phi_S:\Omega\to \R$ as 
\begin{equation*} 
\Phi_S(x):= \int_S \left( \frac1{|x-y|^N}-\frac1{|y|^N} \right)\,\ud y;
\end{equation*}
{in view of 
\eqref{Omega-in-B1},
we have that $\|\Phi_S\|_{L^\infty(\Omega)}\le C(S)$, for some constant $C(S)$ depending only on $S$}. We also set $\Phi_{P,Q}:=\Phi_P-\Phi_Q$.

\begin{prop}
\label{prop_win-hol_bilanciate}
 The set  $G$ defined in \eqref{def_gpq} is balanced, i.e., satisfies \eqref{balcase}
and the limit in~\eqref{lamG} exists with
\begin{equation}\label{lambda_win-hol}
\lambda_G = \int_P\frac{\ud y}{|y|^N} - \int_Q\frac{\ud y}{|y|^N}.
\end{equation}
Moreover, for every $E\in\adm$, setting~$F:=E\cap\Om$, one has
\begin{equation}\label{per0_win-hol}
\Per^0_G(E) = \Per^0_{G_A}(F\cup G_A) + \int_\Om\Phi_{P,Q}\,\ud x - 2\int_F\Phi_{P,Q}\,\ud x.
\end{equation}
Consequently, minimizing~$\Per^0_G$ in $\adm$ is equivalent to minimize
\begin{equation*}
\min\Big\{\Per^0_{G_A}(F\cup G_A)-2\int_F\Phi_{P,Q}\,\ud x\,:\,F\in\M(\R^N),\quad F\subset\Omega\Big\}.
\end{equation*}
\end{prop}

\begin{proof}
{Since $P$ and $Q$ have finite measure, by \eqref{pertufin2}, $\alpha(G)=\alpha(G_A)$, and hence, by \eqref{alfa1cono} $G$ is balanced.}
Furthermore,
$$
\alpha^s_1(G) = \alpha^s_1(G_A)
+ s\int_P\frac{\ud y}{|y|^{N+s}} - s\int_Q\frac{\ud y}{|y|^{N+s}},
$$
whence,  {subtracting~$\alpha(G)$, dividing by $s$ }and letting~$s\to0^+$, we
obtain \eqref{lambda_win-hol} by Dominated Convergence.

It remains to prove \eqref{per0_win-hol}. Since~$P,Q\subset\Cc B_1$,
only the tail part of the functional changes when one passes from~$G_A$ to~
$G$. {More precisely, using also that $G_A\subset\Cc B_1$
\begin{equation*}
\begin{aligned}
\Per^0_G(E)-\Per^0_{G_A}(F\cup G_A)=&\,\phantom{+}\int_{F}\Phi_Q\,\ud x+\int_F\Phi_{(\Cc G_A\setminus P) \cap \Cc B_1}\,\ud x-\int_{F}\Phi_{\Cc G_A \cap \Cc B_1}\,\ud x\\
&\,+\int_{\Cc F\cap\Omega}\Phi_P\,\ud x+\int_{\Cc F\cap\Omega} \Phi_{G_A\setminus Q}\,\ud x-\int_{\Cc F\cap\Omega} \Phi_{G_A}\,\ud x\\
=&\,-\int_{F}\Phi_{P,Q}\,\ud x+\int_{\Cc F\cap\Omega}\Phi_{P,Q}\,\ud x,
\end{aligned}
\end{equation*}
whence \eqref{per0_win-hol} follows.
}
\end{proof}

\begin{rem} 
The construction above also shows that the scalar first-order coefficient
$\lambda_G$ can be varied independently of the leading tail mass. Let
$
G_0:=\{x\in\Cc B_1:\ x_N<0\}$.
Then $G_0$ is balanced and $\lambda_{G_0}=0$.
Denoting by~$d_g$ the geodesic distance on~$\mathbb S^{N-1}$, set
\begin{equation}\label{calotta}
A^\pm_{\beta,\rho}:= \left\{ x\in\R^N:\ 1<|x|<\rho,\ 
d_g\left(\frac{x}{|x|},\pm e_N\right)\le\beta(|x|) \right\},
\end{equation}
where $\rho\in[1,+\infty]$ and
$\beta:(1,+\infty)\to(0,+\infty)$ is measurable. If
$\beta\equiv\bar\beta>0$ and $\rho<+\infty$, then
$$
\int_{A^\pm_{\bar\beta,\rho}}\frac{\ud y}{|y|^N}
=
\mathcal{H}^{N-1}
\left( \left\{\theta\in\mathbb S^{N-1}: d_g(\theta,\pm e_N)\le\bar\beta\right\} \right)
\log\rho .
$$
This quantity is continuous in~$\rho$, vanishes as~$\rho\to1^+$, and tends
to~$+\infty$ as~$\rho\to+\infty$. Hence, by adding a bounded cap contained
in~$\Cc G_0$, or by removing a bounded cap contained in~$G_0$, one obtains
balanced exterior data with arbitrary prescribed finite value of~$\lambda_G$.

The same mechanism also produces the values $\lambda_G=+\infty$ and
$\lambda_G=-\infty$. {In this case}, it is enough to choose an aperture function
$\beta(r)\to0$ as $r\to+\infty$ so slowly that
$$
\int_1^{+\infty} \mathcal{H}^{N-1}
\left( \left\{\theta\in\mathbb S^{N-1}: d_g(\theta,\pm e_N)\le\beta(r)\right\}
\right) \frac{\ud r}{r}=+\infty.
$$
Then $\alpha(A^\pm_{\beta,+\infty})=0$, while the corresponding logarithmic
first-order contribution is infinite.
\end{rem}

\medskip
\subsection{Comparison with the local limit}

We finally compare the tail selection described above with the local limit
$s\to1^-$. The two endpoint regimes detect different parts of the exterior
datum. As $s\to0^+$, the limit is governed by the tail at infinity and, at
the first nontrivial order, by the renormalized geometry encoded in
$\Per^0_G$. As $s\to1^-$, the limit is local and sees only the trace induced
near~$\partial\Om$.

The balanced windows-and-holes construction makes the distinction explicit.
The bounded perturbations $P$ and~$Q$ in
Proposition~\ref{prop_win-hol_bilanciate} do not change the leading tail
mass, but they may select a droplet through the potential~$\Phi_{P,Q}$.
If these perturbations stay at positive distance from~$\Om$, they are
invisible to the local limit. Thus the first-order expansion at $s=0$ and
the local limit at $s=1$ encode genuinely different geometric information.

\medskip
The following examples illustrate this contrast for varying exterior data.
Example \ref{esempio-variabil} concerns the unconstrained problem and shows the dominance of
the tail mass as~$s\to0^+$ whereas
Example \ref{esempio-variabil2} focuses on the volume-constrained
problem and shows the additional boundary-layer dependence of diagonal local
limits as $s\to1^-$.
\begin{example}
\label{esempio-variabil}
Let $\Om=B_1$. Recalling the notation in~\eqref{calotta}, for every
$k\in\N$ we set
$$
G_k:=\Cc B_{1+\frac1k}\cup A^-_{1,+\infty}.
$$
The sets $G_k$ differ from $\Cc B_1$ only in the bounded annulus
$B_{1+\frac1k}\setminus B_1$, outside the cap $A^-_{1,+\infty}$. Hence they
are bounded perturbations of the datum $\Cc B_1$ in the sense of
Remark~\ref{rem:stable}. In particular the leading tail mass is the maximal
one,
$$
\alpha(G_k)=N\omega_N,
$$
at the scale relevant for the zero-order limit. Therefore, by
Remark~\ref{rem:stable}, or - equivalently - by the varying-datum version of
Proposition~\ref{casodivunmezzo}, for every sequence $s_k\to0^+$ and every
sequence $\{E_k\}_{k\in\N}$ of minimizers of
$$
\min\{\Per^{s_k}_{G_k}(E):\ E\in\mathscr A_{G_k}\},
$$
one has
$$
\chi_{E_k\cap B_1}\to1
\qquad\textrm{in }L^1(B_1).
$$
Thus, in the tail regime $s\to0^+$, the exterior datum forces the full phase
inside the domain.

\medskip
The behavior in the local regime is different. Let
$$
\Gamma^-:=\overline{A^-_{1,+\infty}}\cap\partial B_1
$$
be the spherical trace induced on $\partial B_1$ by the cap
$A^-_{1,+\infty}$. By Theorem~\ref{s1} and a standard diagonal argument,
there exists a sequence $\{t_k\}_{k\in\N}$ with $t_k\to1^-$ sufficiently fast,  such that every sequence $\{E_k\}_{k\in\N}$ of minimizers of
$$
\min\{\Per^{t_k}_{G_k}(E):\ E\in\mathscr A_{G_k}\}
$$
converges, in $L^1_{\rm loc}(\R^N)$, to
$$
\Cc B_1\cup K^-,
$$
where $K^-\subset\overline{B_1}$ is the Euclidean perimeter-minimizing hull
with trace $\Gamma^-$ on $\partial B_1$. In the present symmetric situation,
$K^-$ is precisely the convex hull in $\overline{B_1}$ of the spherical cap
$\Gamma^-$.

Thus the same sequence of exterior data produces two different limiting
selections. In the tail regime $s\to0^+$, the maximal tail mass fills the
whole free region $B_1$. In the diagonal local regime $s\to1^-$, the far
tail is invisible and the selection is governed only by the trace induced on
$\partial B_1$.
\end{example}
Before providing the example on the different asymptotic behavior, as $s\to 0^+$ and $s\to 1^-$, of the minimizers for the volume constrained problems, we focus on the extremal cases $G=\emptyset$ and $G=\Cc\Omega$.
\begin{rem}
\label{rem:ovvietaB}
Let $0<s<1$.
If $G=\emptyset$, then 
every admissible set is contained in~$\Om$ up to
a null set, and
$$
\Per^s_G(E)=\I^s(E,\Cc E).
$$
Thus the constrained problem reduces to the standard fractional isoperimetric
problem for sets confined in~$\Om$. If
$$
m_\Om:=\sup\{|B|:\ B\textrm{ is a ball and }B\subset\Om\},
$$
then, for every~$m\in(0,m_\Om]$, the minimizers are the balls of volume
$m$ contained in~$\Om$.

Analogously, if $G=\Cc\Om$, then the same statement applies to the
complementary phase. More precisely, for every~$m\in(0,m_\Om]$, the
minimizers of
\begin{equation}\label{minze}
\min\{\Per^s_{\Cc\Omega}(E):\ E\in\mathscr{A}_{\Cc\Omega},\ |\Cc E\cap\Om|=m\}
\end{equation}
are the sets
\begin{equation}\label{claminze}
\{\Cc B:\ B\textrm{ is a ball},\ B\subset\Om,\ |B|=m\}.
\end{equation}
At the first-order level, Theorem~\ref{teovol0} gives the corresponding
minimizers of $\Per^0_{\emptyset}$ and $\Per^0_{\Cc\Om}$ under the
same volume constraints. In view of Proposition~\ref{comezeroper} and
Remark~\ref{rem:neiminimi}, the uniqueness and stability of these minimizers
are related to the stability result in~\cite[Remark 8.3]{DLNP21} {(see also \cite{C17,FL21})}, {that shows in particular that for $s=0$  the class of minimizers of \eqref{minze} coincides with \eqref{claminze}}.
\end{rem}
\begin{example} \label{esempio-variabil2}
Let $\Om=B_1$ and let $m\in(0,\omega_N)$. For every $k\in\N$, we set
$$
G_k:=\Cc B_{1+\frac1k}.
$$
Let $s_k\to0^+$, and let $\{E_k\}_{k\in\N}$ be a sequence of minimizers of
$$
\min\left\{ \Per^{s_k}_{G_k}(E):\ E\in\mathscr A_{G_k},\ |E\cap B_1|=m \right\}.
$$
By Theorem~\ref{teoval0var}
 and Remark~\ref{rem:ovvietaB}, every sequence of
minimizers converges, up to subsequences, to a set of the form
$$
\Cc B_r(\xi), \qquad \omega_N(1-r^N)=m, \qquad B_r(\xi)\subset B_1.
$$
Thus, in the tail regime~$s\to0^+$, the limiting configuration is selected
as in the complementary confined problem.

The local regime~$s\to1^-$ is sensitive to the competition between the
parameter $1-s$ and the distance of~$G_k$ from~$\partial B_1$. On the one hand,
one can choose a sequence~$\{s_k\}_{k\in\N}$ with $s_k\to1^-$ sufficiently fast, so that the exterior datum is seen as separated from~$\partial B_1$ at the local scale. In this case,  by
Theorem~\ref{s1} and by a standard diagonal argument, the limits (as $k\to +\infty$) of minimizers of $\Per^{s_k}_{G_k}$ in $\mathscr{A}_{G_k,m}$
are union of~$\Cc B_1$ with a ball of measure $m$ contained in~$B_1$. 
On the other hand, choosing a sequence~$\{s_k\}_{k\in\N}$ with $s_k\to1^-$ sufficiently slow,   the limit is governed by the relative
isoperimetric problem in~$B_1$. In such a case, the limiting set inside~$B_1$ is a
minimizer of the Euclidean relative perimeter among sets of volume~$m$;
equivalently, it is a spherical cap meeting~$\partial B_1$ orthogonally {(see \cite{BS79})}.

Therefore the varying-datum problem has different limiting selections in the
tail regime $s\to0^+$ and in diagonal local regimes~$s\to1^-$. This
confirms that the two endpoint limits do not encode the same geometric
information.
\end{example}

\appendix
\section{Proof of Proposition \ref{comezeroper}}\label{appendice}
We conclude this paper by proving the  {connection of the first-order $\Gamma$-limit with the $0$-fractional perimeter $\hat H^0$.}

\begin{proof}[Proof of Proposition \ref{comezeroper}]
We give the proof in detail; the argument is only a decomposition of the renormalized interaction, but the cancellation at infinity has to be kept explicit. {To this end we} separate the interactions according to whether the exterior variable lies in~$B_1$ or in its complement. This gives
\begin{equation}\label{thecla1}
\begin{aligned}
\Per^0_{\emptyset}(E)=&\,\phantom{+}\I^0(E,\Cc E\cap\Om)-|E|\int_{\Om\cap \Cc B_1}\frac{\ud y}{|y|^N}\\
&\,+\I^0(E,\Cc\Om\cap B_1)+\int_{E}\ud x\int_{\Cc \Om\cap\Cc B_1}\Big(\frac{1}{|x-y|^N}-\frac{1}{|y|^N}\Big)\,\ud y.
\end{aligned}
\end{equation}
Moreover, since $\Cc E=(\Cc E\cap \Om)\cup\Cc\Om$ (where the union is disjoint), we have that for every $x\in\R^N$
$$
\Cc E\cap B_1(x)=(\Cc E\cap \Om\cap B_1(x))\cup (\Cc\Om\cap B_1(x)),
$$
whence we deduce
\begin{equation}\label{thecla2}
\begin{aligned}
\I^0(E,\Cc E\cap\Om)
=&\,\phantom{+}\int_{E}\ud x\int_{\Cc E\cap\Om\cap B_1(x)}\frac{\ud y}{|x-y|^N}
+
\int_E\ud x\int_{\Cc E\cap \Om\cap\Cc B_1(x)}\frac{\ud y}{|x-y|^N}\\
=&\,\phantom{+}\int_{E}\ud x\int_{\Cc E\cap B_1(x)}\frac{\ud y}{|x-y|^N}+
\int_E\ud x\int_{\Cc E\cap \Om\cap\Cc B_1(x)}\frac{\ud y}{|x-y|^N}\\
&\,-\int_{E}\ud x\int_{\Cc \Om\cap B_1(x)}\frac{\ud y}{|x-y|^N};
\end{aligned}
\end{equation}
{furthermore,}
\begin{equation}\label{thecla20}
\begin{aligned}
\I^0(E,\Cc\Om\cap B_1)
=&\,\int_{E}\ud x\int_{\Cc\Om\cap B_1\cap B_1(x)}\frac{\ud y}{|x-y|^N}+\int_{E}\ud x\int_{\Cc\Om\cap B_1\cap\Cc B_1(x)}\frac{\ud y}{|x-y|^N}.
\end{aligned}
\end{equation}
{By \eqref{thecla2} and \eqref{thecla20}, we have} 
\begin{equation}\label{thecla200}
\begin{aligned}
&\,\phantom{+}\I^0(E,\Cc E\cap\Om)+\I^0(E,\Cc\Om\cap B_1)\\
=&\,\phantom{+}\int_{E}\ud x\int_{\Cc E\cap B_1(x)}\frac{\ud y}{|x-y|^N}
+
\int_E\ud x\int_{\Cc E\cap \Om\cap\Cc B_1(x)}\frac{\ud y}{|x-y|^N}\\
&\,+\int_{E}\ud x\int_{\Cc\Om\cap B_1\cap\Cc B_1(x)}\frac{\ud y}{|x-y|^N}-\int_{E}\ud x\int_{\Cc \Om\cap B_1(x)\cap \Cc B_1}\frac{\ud y}{|x-y|^N}.
\end{aligned}
\end{equation}
Analogously,
\begin{equation*}\label{thecla21}
\begin{aligned}
-|E|\int_{\Om\cap\Cc B_1}\frac{\ud y}{|y|^N} 
=&\,-\int_{E}\ud x\int_{\Om\cap \Cc B_1\cap \Cc B_1(x)}\frac{\ud y}{|y|^N}
-\int_{E}\ud x\int_{\Om\cap \Cc B_1\cap B_1(x)}\frac{\ud y}{|y|^N}\\
=&\,\phantom{-}\int_E\ud x\int_{\Om\cap \Cc B_1\cap\Cc B_1(x)}
\Big(\frac{1}{|x-y|^N}-\frac{1}{|y|^N}\Big)\ud y\\
&\,-\int_{E}\ud x\int_{\Om\cap\Cc B_1\cap\Cc B_1(x)}\frac{\ud y}{|x-y|^N}-\int_{E}\ud x\int_{\Om\cap\Cc B_1\cap B_1(x)}\frac{\ud y}{|y|^N}\\
=&\,\phantom{-}\int_E\ud x\int_{\Om\cap \Cc B_1\cap\Cc B_1(x)}
\Big(\frac{1}{|x-y|^N}-\frac{1}{|y|^N}\Big)\ud y\\
&\,-\int_{E}\ud x\int_{\Om\cap\Cc B_1(x)}\frac{\ud y}{|x-y|^N}+\int_{E}\ud x\int_{\Om\cap B_1\cap\Cc B_1(x)}\frac{\ud y}{|x-y|^N}\\
&\,-\int_{E}\ud x\int_{\Om\cap\Cc B_1\cap B_1(x)}\frac{\ud y}{|y|^N},
\end{aligned}
\end{equation*}
and
\begin{equation*}\label{thecla210}
\begin{aligned}
\int_{E}\ud x\int_{\Cc \Om\cap\Cc B_1}\Big(\frac{1}{|x-y|^N}-\frac{1}{|y|^N}\Big)\,\ud y
=&\,\phantom{+}\int_{E}\ud x\int_{\Cc \Om\cap\Cc B_1\cap\Cc B_1(x)}\Big(\frac{1}{|x-y|^N}-\frac{1}{|y|^N}\Big)\,\ud y\\
&\,+\int_{E}\ud x\int_{\Cc \Om\cap\Cc B_1\cap B_1(x)}\Big(\frac{1}{|x-y|^N}-\frac{1}{|y|^N}\Big)\,\ud y,
\end{aligned}
\end{equation*}
whence we deduce
\begin{equation}\label{thecla3}
\begin{aligned}
&\,-|E|\int_{\Om\cap\Cc B_1}\frac{\ud y}{|y|^N}+\int_{E}\ud x\int_{\Cc \Om\cap\Cc B_1}\Big(\frac{1}{|x-y|^N}-\frac{1}{|y|^N}\Big)\,\ud y\\
=&\,\phantom{-}\int_{E}\ud x\int_{\Cc B_1\cap\Cc B_1(x)}\Big(\frac{1}{|x-y|^N}-\frac{1}{|y|^N}\Big)\,\ud y\\
&\,-\int_{E}\ud x\int_{\Om\cap\Cc B_1(x)}\frac{\ud y}{|x-y|^N}+\int_{E}\ud x\int_{\Om\cap B_1\cap\Cc B_1(x)}\frac{\ud y}{|x-y|^N}\\
&\,+\int_{E}\ud x\int_{\Cc \Om\cap\Cc B_1\cap B_1(x)}\frac{\ud y}{|x-y|^N}-\int_{E}\ud x\int_{{\Cc B_1\cap B_1(x)}}\frac{\ud y}{|y|^N}\\
=&\,-\int_{E}\ud x\int_{E\cap\Cc B_1(x)}\frac{\ud y}{|x-y|^N}-\int_{E}\ud x\int_{\Cc E\cap\Om\cap\Cc B_1(x)}\frac{\ud y}{|x-y|^N}\\
&\,+\int_{E}\ud x\int_{\Om\cap B_1\cap\Cc B_1(x)}\frac{\ud y}{|x-y|^N}+\int_{E}\ud x\int_{\Cc \Om\cap\Cc B_1\cap B_1(x)}\frac{\ud y}{|x-y|^N}\\
&\,-\int_{E}\ud x\int_{\Cc B_1\cap B_1(x)}\frac{\ud y}{|y|^N},
\end{aligned}
\end{equation}
where in the last equality we used the change of variable $z=x-y$ to deduce that
\begin{equation*}
\begin{aligned}
&\,\int_{E}\ud x\int_{\Cc B_1\cap\Cc B_1(x)}\Big(\frac{1}{|x-y|^{N}}-\frac{1}{|y|^N}\Big)\,\ud y\\
=&\,\int_{E}\ud x\int_{\Cc B_1(x)\cap\Cc B_1}\Big(\frac{1}{|z|^{N}}-\frac{1}{|x-z|^N}\Big)\,\ud z
=0.
\end{aligned}
\end{equation*}
By \eqref{thecla1}, \eqref{thecla200}, \eqref{thecla3} and \eqref{defH0}, we obtain
\begin{equation*}
\begin{aligned}
\Per^0_{\emptyset}(E)=&\,\phantom{+}\hat{H}^0(E)+\int_{E}\ud x\int_{B_1\cap \Cc B_1(x)}\frac{\ud y}{|x-y|^N}-\int_{E}\ud x\int_{\Cc B_1\cap B_1(x)}\frac{\ud y}{|y|^N}\\
=&\,\phantom{+}\hat{H}^0(E),
\end{aligned}
\end{equation*}
where in the last equality we have used  the change of variable $z=x-y$
to  deduce that
\begin{equation*}
\int_{E}\ud x\int_{\Cc B_1\cap B_1(x)}\frac{\ud y}{|y|^N}
=\int_{E}\ud x\int_{\Cc B_1(x)\cap B_1}\frac{\ud z}{|z-x|^N}.
\end{equation*}
\end{proof}

\bigskip

\end{document}